\documentclass{article}
\usepackage[english]{babel}
\usepackage{amsmath,amssymb,eufrak,amsthm}
\begin{document}
\bibliographystyle{plain}

\title{Ideals of varieties parameterized by certain symmetric tensors}
\author{Alessandra Bernardi}
\date{}
\maketitle

\renewcommand{\sectionmark}[1]{\markboth{#1}{}}

\def\move-in{\parshape=1.75true in 5true in}

\def\CC{\mathbb C}
\def\ZZ{\mathbb Z}
\def\NN{\mathbb N}
\def\RR{\mathbb R^}
\def\PP{\mathbb P^}
\def\CA{\mathcal A}
\def\G{\mathfrak{S}}

\def\sgn{\mathrm{sgn}}
\def\dim{\mathrm{dim}}
\def\codim{\mathrm{codim}}
\def\rk{\mathrm{rk}}
 
\def\v{\underline{v}}
\def\w{\underline{w}}
\def\u{\underline{u}}

\def\Def#1{\noindent {\bf Definition #1:}}
\def\Not{\noindent {\bf Notation: }}
\def\Prop#1{\noindent {\bf Proposition #1:}}
\def\Proof{\noindent {\it Proof: }}
\def\Obs{\noindent{\bf Remark: }}
\def\Ex{\noindent{\bf Example: }}
\newenvironment{dem}{\begin{proof}[Proof]}{\end{proof}}
\newtheorem{theorem}{Theorem}[section]
\newtheorem{lemma}[theorem]{Lemma}
\newtheorem{propos}[theorem]{Proposition}
\newtheorem{corol}[theorem]{Corollary}
\newtheorem{defi}[theorem]{Definition}
\newtheorem{conj}[theorem]{Conjecture}
\newtheorem{es}[theorem]{Esercizio}

{\small \begin{quote} ABSTRACT. The ideal of a Segre variety $\PP {n_{1}}\times \cdots \times \PP {n_{t}}\hookrightarrow \PP {(n_{1}+1)\cdots (n_{t}+1)-1}$ is generated by the $2$-minors of a generic hypermatrix of indeterminates (see \cite{Ha1} and \cite{Gr}). We extend this result to the case of Segre-Veronese varieties. The main tool is the concept of ``weak generic hypermatrix'' which allows us to treat also the case of projection of Veronese surfaces from a set of general points and of Veronese varieties  from a Cohen-Macaulay subvariety of codimension $2$.
\end{quote}}

\section{Introduction}

In this paper we study the generators of the ideal of Segre-Veronese varieties and the ideal of projections of Veronese surfaces from a set of general points and, more generally, of Veronese varieties  from a Cohen-Macaulay subvariety of codimension $2$.

A Segre variety parameterizes completely decomposable tensors (Definition \ref{decompos}).

The problem of tensor decomposition has been studied studied for many years and by researchers in many scientific areas as Algebraic Geometry (see for example \cite{CGG1}, \cite{LM}, \cite{LW}, \cite{AOP}, \cite{Za}), Algebraic Statistic (see \cite{HR}, \cite{GSS}, \cite{PS}), Phylogenetic (\cite{AR}, \cite{Bo},  \cite{Lak}, \cite{SS}),  Telecommunications (\cite{Com}), Complexity Theory (\cite{BCS}, \cite{Lan}, \cite{Li}, \cite{St}), Quantum Computing (\cite{BZ}), Psychometrics (\cite{CKP}), Chemometrics (\cite{Br}).

 In \cite{Ha1} (Theorem 1.5) it is proved that the ideal of a Segre variety is generated by all the $2$-minors of a generic hypermatrix of indeterminates. 

Here we prove an analogous statement for Segre-Veronese varieties (see \cite{CGG2}). Segre-Veronese varieties parameterize certain symmetric decomposable tensors, and are the embedding of $\PP {n_{1}}\times \cdots \times \PP {n_{t}}$ into $\PP {\Pi_{i=1}^{t}{n_{i}+d_{i}-1 \choose d_{i}}-1}$ given by the sections of the sheaf ${\mathcal O}(d_{1}, \ldots , d_{t})$ with $d_{1}, \ldots , d_{t}\in \NN$ (see Section \ref{SegreVeronese}). We prove (in Theorem \ref{theor segre ver}) that their ideal is generated by the $2$-minors of a generic symmetric hypermatrix (Definition \ref{hyper sym}).

The idea we use is the following; generalizing ideas in \cite{Ha1} we define ``weak generic hypermatrices'' (see Definition \ref{weak hypermatrix}) and we prove that the ideal generated by $2$-minors of a weak generic hypermatrix is a prime ideal (Proposition \ref{gen1.12}). Then we show that a symmetric hypermatrix of indeterminates is  weak generic and we can conclude, since the ideal generated by its $2$-minors defines, set-theoretically, a Segre-Veronese variety.

An analogous idea is used in Sections \ref{PoVs} and \ref{proiezioni veronese variet}  in order to find the generators of  projections of Veronese varieties from a subvariety of codimension $2$. This is a problem which has been studied classically in Algebraic Geometry (starting with the projection of Veronese surface, see \cite{Sh}); for a quite general analysis of subalgebras of the Rees Algebra associated to embeddings of blow ups of $\PP n$ along subvarieties, see \cite{CHTV} and \cite{MU}.

 Denote with $Y_{n,d}$ the Veronese variety obtained as the $d$-uple embedding of $\PP n$ into $\PP {{n+d\choose d}-1}$ and consider the surface $Y\subset \PP {{2+d \choose 2}-s-1}$ which is the projection of $Y_{2,d}$ from $s$ general points on it. The defining ideal of $Y$ has been studied in \cite{Ha1} when $s$ is a binomial and $s\leq {d\choose 2}$ and in \cite{GL} and \cite{Ha2} for $s>{d\choose 2}$ (in the second paper also the case of any set of $s$ points is treated, when $d\geq \max \{4,s+1\}$). Here we complete the picture for $s< {d\choose 2}$ general points on $Y_{2,d}$; our method follows the framework of \cite{GG} and \cite{GL}, but uses the ``hypermatrix'' point of view of \cite{Ha1}. We construct a hypermatrix in such a way that its $2$-minors together with some linear equations generate an ideal $I$ that defines $Y$ set-theoretically; then we prove that such hypermatrix is weak generic and in Theorem \ref{end surf} we prove that $I$ is actually the ideal of  the projected surface.

This construction can be generalized to projections of Veronese varieties $Y_{n,d}$, for all $n,d>0$, from a subvariety of codimension $2$ and of degree $s={t+1 \choose 2}+k \leq {d\choose 2}$ for some non negative integers $t$, $k$, $d$ such that $0<t<d-1$ and $0\leq k \leq t$ (see Section \ref{proiezioni veronese variet}).
\\
\\
I want to thank A. Gimigliano for many useful talks and suggestions.
\\
I want to thank also J. M. Landsberg for pointing out to me that by
using representation theory techniques (as for example in \cite{LM}) it
is possible to see that the equations coming from the vanishing of
$2$-minors of a symmetric hypermatrix of indeterminates are the
generators of the ideal of a Segre-Veronese variety.
By an unpublished theorem of Kostant, the ideal of any homogeneously embedded
rational homogeneous variety is generated in degree two by the
annihilator of a certain vector space (for the experts:  the homogeneously embedded
rational homogeneous variety $G/P\subset \Bbb P (V_{\lambda})$, is
generated in degree two by the
annihilator of $V_{2\lambda}$ in $S^2(V^*_{\lambda})$). While the
representation-theoretic techniques identify the modules generating the
ideal, they do not provide an explicit method for writing
down a set of generators, which is the subject of this paper.  
\\
Last but not least thanks to the anonymous referee for his careful work and suggestions.

\section{Preliminaries}\label{Prelim}

Let $K=\overline{K}$ be an algebraically closed field of characteristic zero, and let $V_{1}, \ldots , V_{t}$ be vector spaces over $K$ of dimensions $n_{1}, \ldots , n_{t}$ respectively. We will call en element $T \in V_{1}\otimes \cdots \otimes V_{t}$ a tensor of size $n_{1}\times \cdots \times n_{t}$.

Let $E_{j}=\{\underline{e}_{j,1}, \ldots , \underline{e}_{j,n_{j}}\}$ be a basis for the vector space $V_{j}$, $j=1, \ldots , t$. We define a basis $E$ for $V_{1}\otimes \cdots \otimes V_{t}$ as follows:
\begin{equation}\label{canbas}E:=\{\underline{e}_{i_{1}, \ldots , i_{t}}=\underline{e}_{1,i_{1}}\otimes \cdots \otimes \underline{e}_{t,i_{t}} \; | \; 1\leq i_{j} \leq n_{j}, \, \forall j=1, \ldots , t\}.\end{equation}
A tensor $T\in V_{1}\otimes \cdots \otimes V_{t}$ 
can be represented via a so called ``hypermatrix'' (or ``array'') 
$$\CA = (a_{i_{1}, \ldots , i_{t}})_{1\leq i_{j}\leq n_{j}\, , \, j=1, \ldots , t} $$
with respect to the basis $E$ defined in (\ref{canbas}), i.e.:
$$T=\sum_{1\leq i_{j} \leq n_{j}, \, j=1, \ldots , t} a_{i_{1}, \ldots , i_{t}}\underline{e}_{i_{1}, \ldots , i_{t}}.$$

\begin{defi}\label{decompos} A tensor $T\in V_{1}\otimes \cdots \otimes V_{t}$ is called ``decomposable'' if, for all $j=1, \ldots , t$, there exist $\v_{j}\in V_{j}$ such that $T=\v_{1}\otimes \cdots \otimes \v_{t}$.
\end{defi}

\begin{defi} Let $E_{j}=\{\underline{e}_{j,1}, \ldots ,\underline{e}_{j,n_{j}} \}$ be a basis for the vector space $V_{j}$ for $j=1, \ldots , t$. Let also $\underline{v}_{j}=\sum_{i=1}^{n_{j}} a_{j,i}\underline{e}_{j,i}\in V_{j}$ for $j=1, \ldots , t$. The image of the following embedding
$$\begin{array}{rcl}
\PP {}(V_{1})\times \cdots \times \PP {}(V_{t})&\hookrightarrow & \PP {}(V_{1}\otimes \cdots \otimes V_{t})\\
([\v_{1}]\; \; , \; \; \ldots \; \; , \; \; [\v_{t}]) &\mapsto &[\v_{1}\otimes\cdots \otimes \v_{t}]=\\
&&=\sum_{1\leq i_{j}\leq n_{j}, \, j=1, \ldots , t} [(a_{1,i_{1}}\cdots a_{t,i_{t}}) \underline{e}_{i_{1},\ldots ,i_{t}}]
\end{array}$$
is well defined and it is known as ``Segre Variety''. We denote it by $Seg(V_{1}\otimes \cdots \otimes V_{t})$.
 \end{defi}

\Obs A Segre variety $Seg(V_{1}\otimes \cdots \otimes V_{t})$ parameterizes the decomposable tensors of $V_{1}\otimes \cdots \otimes V_{t}$.
\\
\\
A set of equations defining $Seg(V_{1}\otimes \cdots \otimes V_{t})$ is well known (one of the first reference for a set-theoretical description of the equations of Segre varieties is \cite{Gr}). Before introducing that result we need the notion of $d$-minor of a hypermatrix.
\\
\\
\Not
\begin{itemize}
\item The hypermatrix $\CA =(x_{i_{1}, \ldots , i_{t}})_{1\leq i_{j}\leq n_{j}\, , \, j=1, \ldots , t}$ is said to be a \emph{generic hypermatrix of indeterminates} (or more simply \emph{generic hypermatrix}) of $S:=K[x_{i_{1}, \ldots , i_{t}}]_{1\leq i_{j}\leq n_{j}\, , \, j=1, \ldots , t}$, 
 if the entries of $\CA$ are the independent variables of $S$.
\item We denote by $S_{t}$ the homogeneous degree $t$ part of the polynomial ring $S$.
\item We will always suppose that we have fixed a basis $E_{i}$ for each $V_{i}$ and the basis $E$ for $V_{1}\otimes \cdots \otimes V_{t}$ as in (\ref{canbas}).
\item When we will write ``$\CA$ is the hypermatrix associated to the tensor $T$'' (or vice versa) we will always assume that the association is via the fixed basis $E$. Moreover  if the size of $T$ is $n_{1}\times \cdots \times n_{t}$, then $\CA$ is of the same size.
\end{itemize}

It is possible to extend the notion of ``$d$-minor of a matrix'' to that of  ``$d$-minor of a hypermatrix''.

\begin{defi}  Let $V_{1}, \ldots , V_{t}$ be vector spaces of dimensions $n_{1}, \ldots , n_{t}$, respectively, and let $(J_{1},J_{2})$ be a partition of the set $\{1, \ldots , t\}$. If $J_{1}=\{h_{1}, \ldots , h_{s}\}$ and $J_{2}=\{1, \ldots ,t\} \backslash J_{1}=\{k_{1}, \ldots , k_{t-s}\}$, the $(J_{1},J_{2})$-Flattening of $V_{1}\otimes \cdots \otimes V_{t}$  is the following:
$$V_{J_{1}}\otimes V_{J_{2}}=(V_{h_{1}}\otimes \cdots \otimes V_{h_{s}})\otimes (V_{k_{1}}\otimes \cdots \otimes V_{k_{t-s}}).$$
\end{defi}

\begin{defi} Let $V_{J_{1}}\otimes V_{J_{2}}$ be any flattening of $V_{1}\otimes \cdots \otimes V_{t}$ and let  $f_{J_{1},J_{2}}:\PP {}(V_{1}\otimes \cdots \otimes V_{t})\stackrel{\sim}{\rightarrow} \PP {}(V_{J_{1}}\otimes V_{J_{2}})$ be the obvious isomorphism. Let $\CA$ be a hypermatrix associated to a tensor $T\in V_{1}\otimes \cdots \otimes V_{t}$; let $[T']=f_{J_{1},J_{2}}([T])\in \PP {}(V_{J_{1}}\otimes V_{J_{2}})$ and let $A_{J_{1},J_{2}}$ be the matrix associated to $T'$. Then the $d$-minors of the matrix $A_{J_{1},J_{2}}$ are said to be ``$d$-minors of $\CA$''.
\end{defi}

Sometimes we will improperly write ``a $d$-minor of a tensor $T$'', meaning that it is a $d$-minor of the hypermatrix associated to such a tensor via the fixed basis $E$ of $V_{1}\otimes \cdots \otimes V_{t}$.
\\
\\
\Ex \textbf{$d$-minors of a decomposable tensor.}

Let $V_{1}, \ldots , V_{t}$ and  $(J_{1},J_{2})=(\{h_{1}, \ldots , h_{s}\}, \{k_{1}, \ldots , k_{t-s}\})$ as before.
Consider the following composition of maps:
$$\PP {}(V_{1})\times \cdots \times \PP {}(V_{t}) \stackrel{s_{1}\times s_{2}}{\rightarrow} \PP {}(V_{J_{1}})\times \PP {}(V_{J_{2}})\stackrel{s}{\rightarrow} \PP {}(V_{J_{1}}\otimes V_{J_{2}})$$
where  $Im(s_{1}\times s_{2})=Seg(V_{J_{1}})\times Seg(V_{J_{2}})$ and $Im(s)$ is the Segre variety of two factors.

Consider the basis (made as $E$ above) $E_{J_{1}}$  for $V_{J_{1}}$ and $E_{J_{2}}$  for $V_{J_{2}}$. 
In terms of coordinates, the composition $s\circ (s_{1}\times s_{2})$ is described as follows.

 Let $\v_{i}=(a_{i,1}, \ldots , a_{i,n_{i}})\in V_{i}$ for each $i=1, \ldots , t$ and $T=\v_{1}\otimes \cdots \otimes \v_{t}\in V_{1}\otimes \cdots \otimes V_{t}$; then:
$$s_{1}\times s_{2}([(a_{1,1}, \ldots , a_{1,n_{1}})], \ldots , [(a_{t,1}, \ldots , a_{t,n_{t}})])=([(y_{1,\ldots , 1}, \ldots , y_{n_{h_{1}},\ldots , n_{h_{s}}})],[(z_{1, \ldots , 1}, \ldots , z_{n_{k_{1}}, \ldots , n_{k_{t-s}}})])$$
where $y_{l_{1}, \ldots , l_{s}}=a_{h_{1},l_{1}}\cdots a_{h_{s},l_{s}}$, for $l_{m}=1, \ldots , n_{m}$ and $m=1, \ldots , s$; 
\\and $z_{l_{1}, \ldots , l_{t-s}}=a_{k_{1},l_{1}}\cdots a_{k_{t-s},l_{t-s}}$ for $l_{m}=1, \ldots ,n_{m}$ and $m=1, \ldots , t-s$.

 If we rename the variables in $V_{J_{1}}$ and in $V_{J_{2}}$ as: $(y_{1,\ldots , 1}, \ldots , y_{n_{h_{1}},\ldots , n_{h_{s}}})=(y_{1}, \ldots , y_{N_{1}})$, with $N_{1}=n_{h_{1}}\cdots n_{h_{s}}$, and $(z_{1, \ldots , 1}, \ldots , z_{n_{k_{1}}, \ldots , n_{k_{t-s}}})=(z_{1}, \ldots , z_{N_{2}})$, with $N_{2}=n_{k_{1}}\cdots n_{k_{t-s}}$, then:
$$s([(y_{1}, \ldots , y_{N_{1}})], [(z_{1}, \ldots , z_{N_{2}})])=[(q_{1,1},q_{1,2}, \ldots ,q_{N_{1},N_{2}} )]=
s\circ (s_{1}\times s_{2})([T]),$$
where $q_{i,j}=y_{i}z_{j}$ for $i=1, \ldots , N_{1}$ and $j=1, \ldots , N_{2}$. We can easily rearrange coordinates and write $s\circ (s_{1}\times s_{2})([T])$ as a matrix:
\begin{equation}\label{s1s2sT}
((s_{1}\times s_{2})\circ s)([T])=
\left(\begin{array}{ccc}
q_{1,1} &\cdots & q_{1,N_{2}}
\\
\vdots && \vdots \\
q_{N_{1},1}&\cdots &q_{N_{1}, N_{2}}
\end{array}\right) .
\end{equation}

A $d$-minor of the matrix $s\circ (s_{1}\times s_{2})([T])$ defined in (\ref{s1s2sT}) is called a $d$-minor of the  tensor $T$.
\\
\\
\Ex The $2$-minors of a hypermatrix $\CA=(a_{i_{1}, \ldots , i_{t}})_{1\leq i_{j}\leq n_{j}, \, j=1, \ldots , t}$ are all of the form:
$$a_{i_{1}, \ldots , i_{m}, \ldots , i_{t}}a_{l_{1}, \ldots , l_{m}, \ldots , l_{t}}-a_{i_{1}, \ldots , l_{m}, \ldots , i_{t}}a_{l_{1}, \ldots , i_{m}, \ldots , l_{t}}$$ 
for $1\leq i_{j},l_{j}\leq n_{j}, \, j=1, \ldots , t$ and $1\leq m \leq t$.

\begin{defi} Let  $\CA$ be a  hypermatrix  whose entries are in $K[u_{1}, \ldots , u_{r}]$. The ideal $I_{d}(\CA)$ is the ideal generated by all $d$-minors of $\CA$.
\end{defi}

\Ex The ideal of the $2$-minors of a generic hypermatrix $\CA=(x_{i_{1}, \ldots , i_{t}})_{1\leq i_{j}\leq n_{j},\, j=1, \ldots , t}$ is
$$I_{2}(\CA):= (x_{i_{1}, \ldots , i_{l}, \ldots , i_{t}}x_{j_{1}, \ldots , j_{l}, \ldots , j_{t}}-x_{i_{1}, \ldots , j_{l}, \ldots , i_{t}}x_{j_{1}, \ldots , i_{l}, \ldots , j_{t}})_{l=1, \ldots , t; \; 1\leq i_{k}, j_{k}\leq n_{j}\, , \, k=1, \ldots , t}.$$
\\
\\
It is a classical result (see \cite{Gr}) that a set of equations for a  Segre Variety is given by all the $2$-minors of a generic hypermatrix. In fact, as previously obseved, a Segre variety parameterizes decomposable tensors, i.e. all the ``rank one'' tensors.

In \cite{Ha1} (Theorem 1.5) it is proved that, if $\CA$ is a generic hypermatrix of a polynomial ring $S$ of size $n_{1}\times \cdots \times n_{t}$, then $I_{2}(\CA)$ is a prime ideal in $S$, therefore:
$$I(Seg(V_{1}\otimes \cdots \otimes V_{t}))=I_{2}(\CA)\subset S.$$

Now we generalize this result to another class of decomposable tensors: those defining ``Segre-Veronese varieties''.

\section{Segre-Veronese varieties}\label{SegreVeronese}

\subsection{Definitions and Remarks}

Before defining a Segre-Veronese variety we recall that a Veronese variety  $Y_{n,d}$   is the $d$-uple embedding of $\PP {n}$ into $\PP {{n+d\choose d}-1}$, via  the linear system associated to the sheaf ${\mathcal{O}}(d)$,  with $d>0$.

\begin{defi} A hypermatrix $\CA=(a_{i_{1}, \ldots , i_{d}})_{1\leq i_{j}\leq n, \, j=1, \ldots , d}$ is said to be ``supersymmetric'' if $a_{i_{1}, \ldots , i_{d}}=a_{i_{\sigma(1)}, \ldots , i_{\sigma(d)}}$ for all $\sigma\in \G_{d}$ where $\G_{d}$ is the permutation group of $\{1, \ldots ,d\}$.
\end{defi}

With an abuse of notation we will say that a tensor $T\in V^{\otimes d}$ is supersymmetric if it can be represented by a supersymmetric hypermatrix.

\begin{defi} Let $H\subset V^{\otimes d}$ be the ${n+d-1\choose d}$-dimensional subspace of the supersymmetric tensors of $V^{\otimes d}$, i.e. $H$ is isomorphic to the symmetric algebra $Sym_{d}(V)$.  Let $\tilde{S}$ be a ring of coordinates on $\PP {{n+d-1\choose d}-1}=\PP {}(H)$ obtained as the quotient $\tilde S=S/I$ where $S=K[x_{i_{1}, \ldots , i_{d}}]_{1\leq i_{j}\leq n, \, j=1, \ldots ,d}$  and $I$ is the ideal generated by all  $x_{i_{1}, \ldots , i_{d}}-x_{i_{\sigma(1)}, \ldots , i_{\sigma(d)}}, \forall\; \sigma \in \G_{d}$.
\\
The hypermatrix $(\overline{x}_{i_{1}, \ldots , i_{d}})_{1\leq i_{j}\leq n, \, j=1, \ldots , d}$ whose entries are the indeterminates of $\tilde S$, is said to be a ``generic supersymmetric hypermatrix''.
\end{defi}

\Obs  The Veronese variety $Y_{n-1,d}\subset \PP {{n+d-1\choose d}-1}$ can be viewed as $Seg(V^{\otimes d})\cap \PP {}(H) \subset \PP {}(H)$.\\
Let $\CA=(x_{i_{1}, \ldots  , i_{d}})_{1\leq i_{j}\leq n, \, j=1, \ldots , d} $ be a generic supersymmetric hypermatrix, then it is a known result   that:
\begin{equation}\label{veronese}I(Y_{n-1,d})=I_{2}(\CA)\subset \tilde{S}.\end{equation}
See \cite{Wa} for set theoretical point of view. In \cite{Pu} the author proved that $I(Y_{n-1,d})$ is generated by the $2$-minors of a particular catalecticant matrix (for a definition of ``Catalecticant matrices'' see e.g.  either \cite{Pu} or \cite{Ge}). A. Parolin, in his PhD thesis (\cite{Pa}), proved that the ideal generated by the $2$-minors of that catalecticant matrix is actually $I_{2}(\CA)$, where $\CA$ is a generic supersymmetric hypermatrix. 
\\
\\
In this way we have recalled two very related  facts:
\begin{itemize}
\item if $\CA$ is a generic $n_{1}\times \cdots \times n_{t}$ hypermatrix, then the ideal of the $2$-minors of $\CA$ is the ideal of the Segre variety $Seg(V_{1}\otimes \cdots \otimes V_{t})$;
\item  if $\CA$ is a generic supersymmetric  $\underbrace{n\times \cdots \times n}_{d}$ hypermatrix, then the ideal of the $2$-minors of $\CA$ is the ideal of the Veronese variety $Y_{n-1,d}$, with $\dim (V)=n$.
\end{itemize}

Now we want to prove that a similar result holds also for other kinds of hypermatrices strictly related with those representing tensors parameterized by Segre varieties and Veronese varieties.

\begin{defi} Let $V_{1}, \ldots , V_{t}$ be vector spaces of dimensions $n_{1}, \ldots , n_{t}$ respectively. The Segre-Veronese variety ${\cal S}_{d_{1},\ldots , d_{t}}(V_{1}\otimes \cdots \otimes V_{t})$ is the embedding of    $\PP {}(V_{1})\otimes \cdots \otimes \PP {}(V_{t})$ into $\PP {N-1}$, where $N=\left(\Pi_{i=1}^{t}{n_{i}+d_{i}-1\choose d_{i}}\right)$, given by sections of the sheaf ${\mathcal{O}}(d_{1}, \ldots , d_{t})$. 
\\
I.e. ${\cal S}_{d_{1},\ldots , d_{t}}(V_{1}\otimes \cdots \otimes V_{t})$ is the image of the composition of the following two maps:
$$
\PP {}(V_{1})\times \cdots \times \PP {}(V_{t}) \stackrel{\nu_{d_{1}}\times \cdots \times \nu_{d_{t}}}{\longrightarrow} \PP {{n_{1}+d_{1}-1\choose d_{1}}-1}\times \cdots \times \PP {{n_{t}+d_{t}-1\choose d_{t}}-1}$$
and
$$ \PP {{n_{1}+d_{1}-1\choose d_{1}}-1}\times \cdots \times \PP {{n_{t}+d_{t}-1\choose d_{t}}-1}\stackrel{s}{\longrightarrow} \PP {N-1}$$
where $Im (\nu_{1}\times \cdots \times \nu_{t})=Y_{n_{1}-1,d_{1}}\times \cdots \times Y_{n_{t}-1,d_{t}}$ and $Im (s)$ is the Segre variety with $t$ factors.
\end{defi}

\Ex If $(d_{1}, \ldots ,d_{t})=(1, \ldots , 1)$ then ${\cal S}_{1,\ldots , 1}(V_{1}\otimes \cdots \otimes V_{t})= Seg(V_{1}\otimes \cdots \otimes V_{t})$.
\\
\\
\Ex If $t=1$ and $\dim(V)=n$, then ${\cal S}_{d}(V)$ is the Veronese variety $Y_{n-1,d}$.
\\
\\
Below we describe how to associate to each element of  ${\cal S}_{d_{1}, \ldots , d_{t}}(V_{1}\otimes \cdots \otimes V_{t})$ a decomposable tensor $T\in V_{1}^{\otimes d_{1}}\otimes \cdots \otimes V_{t}^{\otimes d_{t}}$.

\begin{defi} Let $\underline{n}=(n_{1}, \ldots , n_{t})$ and $\underline{d}=(d_{1},\ldots , d_{t})$. If $V_{i}$ are vector spaces of dimension $n_{i}$ for $i=1, \ldots , t$, an ``$(\underline{n},\underline{d})$-tensor'' is defined to be a tensor $T$ belonging to $V_{1}^{\otimes d_{1}}\otimes \cdots \otimes V_{t}^{\otimes d_{t}}$.
\end{defi}

\begin{defi}\label{hyper sym} Let $\underline{n}$  
 and $\underline{d}$ 
 as above. A hypermatrix $\CA=(a_{i_{1,1}, \ldots , i_{1,d_{1}}; \ldots ; i_{t,1}, \ldots , i_{t,d_{t}}})_{1\leq i_{j,k}\leq n_{j}, \, k=1, \ldots , d_{j}, \, j=1, \ldots , t}$ is said to be ``$(\underline{n}, \underline{d})$-symmetric'' if  $a_{i_{1,1}, \ldots , i_{1,d_{1}}; \ldots ; i_{t,1}, \ldots , i_{t,d_{t}}}=a_{i_{\sigma_{1}(1,1)}, \ldots , i_{\sigma_{1}(1,d_{1})}; \ldots ; i_{\sigma_{t}(t,1)}, \ldots , i_{\sigma_{t}(t,d_{t})}}$ for all permutations $\sigma_{j}\in \G (j,d_{j})$ where $\G (j,d_{j})\simeq \G_{d_{j}}$ is the permutation group on $\{(j,1), \ldots , (j,d_{j})\}$ for all $j=1, \ldots , t$.
\end{defi}

An $(\underline{n}, \underline{d})$-tensor $T\in V_{1}^{\otimes d_{1}}\otimes \cdots \otimes V_{t}^{\otimes d_{t}}$ is said to be an ``$(\underline{n}, \underline{d})$-symmetric tensor'' if it can be represented by an $(\underline{n}, \underline{d})$-symmetric hypermatrix.

\begin{defi}\label{Sd1dt}Let $H_{i}\subset V_{i}^{\otimes d_{i}}$ be the subspace of supersymmetric tensors of $V_{i}^{\otimes d_{i}}$ for each $i=1, \ldots , t$, then $H_{1}\otimes \cdots \otimes H_{t}\subset V_{1}^{\otimes d_{1}}\otimes \cdots \otimes V_{t}^{\otimes d_{t}}$ is the subspace of the $(\underline{n}, \underline{d})$-symmetric tensors of $V_{1}^{\otimes d_{1}}\otimes \cdots \otimes V_{t}^{\otimes d_{t}}$. Let $\underline{n}=(n_{1}, \ldots , n_{t})$ and $\underline{d}=(d_{1},\ldots , d_{t})$ and let $R_{[\underline{n}, \underline{d}]}$
be the ring of coordinates on $\PP {N-1}=\PP {}(H_{1}\otimes \cdots \otimes H_{t})$, with $N=\left( \Pi_{i=1}^{t} {n_{i}+d_{i}-1\choose d_{i}}\right)$, obtained from $S=K[x_{i_{1,1}, \ldots , i_{1,d_{1}}; \ldots ; i_{t,1}, \ldots , i_{t,d_{t}}}]_{1\leq i_{j,k}\leq n_{j}, \, k=1, \ldots , d_{j}, \, j=1, \ldots , t}$ via the quotient modulo  $x_{i_{1,1}, \ldots , i_{1,d_{1}}; \ldots ; i_{t,1}, \ldots , i_{t,d_{t}}}-x_{i_{\sigma_{1}(1,1)}, \ldots , i_{\sigma_{1}(1,d_{1})}; \ldots ; i_{\sigma_{t}(t,1)}, \ldots , i_{\sigma_{t}(t,d_{t})}}$, for all $\sigma_{j}\in \G (j,d_{j})$ and $j=1, \ldots , t$.
\\
The hypermatrix $(\overline{x}_{i_{1,1}, \ldots , i_{1,d_{1}}; \ldots ; i_{t,1}, \ldots , i_{t,d_{t}}})_{1\leq i_{j,k}\leq n_{j}, \, k=1, \ldots , d_{j}, \, j=1, \ldots , t}$ of indeterminates of $R_{[\underline{n},\underline{d}]}$, is said to be a ``generic $(\underline{n}, \underline{d})$-symmetric hypermatrix''.
\end{defi}

\Obs  
It is not difficult to check that, as sets: 
\begin{equation}\label{SegVerSet}\PP {}(H_{1}\otimes \cdots \otimes H_{t})\cap Seg(V_{1}^{\otimes d_{1}}\otimes \cdots \otimes V_{t}^{\otimes d_{t}})= {\cal S}_{d_{1}, \ldots , d_{t}}(V_{1}\otimes \cdots \otimes V_{t});\end{equation}
i.e.  ${\cal S}_{d_{1}, \ldots , d_{t}}(V_{1}\otimes \cdots \otimes V_{t})$ parameterizes the  $(\underline{n},\underline{d})$-symmetric decomposable $(\underline{n},\underline{d})$-tensors of $V_{1}^{\otimes d_{1}}\otimes \cdots \otimes V_{t}^{\otimes d_{t}}$.
\\
Since Segre variety is given by the vanishing of $2$-minors of a hypermatrix of indeterminates and $H_{1}\otimes \cdots \otimes H_{i}$ is a linear subspace of $V_{1}\otimes \cdots \otimes V_{t}$, it follows that a Segre-Veronese variety is set-theoretically given by the $2$-minors of an $(\underline{n},\underline{d})$-symmetric hypermatrix of indeterminates .
\\
\\
In Section \ref{Ideal Segre veronese} we will prove that the ideal of the $2$-minors of the generic $(\underline{n},\underline{d})$-symmetric hypermatrix in $R_{[\underline{n}, \underline{d}]}$ is the ideal of a Segre-Veronese variety.
We will need the notion of  ``weak generic hypermatrices'' that we are going to introduce.

\subsection{Weak Generic Hypermatrices}

The aim of this section is Proposition \ref{gen1.12} which asserts that the ideal generated by $2$-minors of a weak generic hypermatrix (Definition \ref{weak hypermatrix}) is prime.

\begin{defi}
A $k$-th section of a hypermatrix $\CA=(x_{i_{1}, \ldots , i_{t}})_{1\leq i_{j}\leq n_{j},  j=1, \ldots , t}$ is a hypermatrix of the form
$$\CA_{i_{k}}^{(l)}=(x_{i_{1}, \ldots , i_{t}})_{1\leq i_{j}\leq n_{j}, j=1,\ldots , \hat{k}, \ldots , t, i_{k}=l}.$$
\end{defi}

\Obs If a hypermatrix $\CA$ represents a tensor $T\in V_{1}\otimes \cdots \otimes V_{t}$, then a $k$-th section of $\CA$ is a hypermatrix representing a tensor $T'\in V_{1}\otimes \cdots \otimes \hat{V}_{k}\otimes \cdots \otimes V_{t}$.

We introduce now the notion of ``weak generic hypermatrices''; this is a generalization of ``weak generic box'' in \cite{Ha1}.

\begin{defi}\label{weak hypermatrix} Let $K[u_{1}, \ldots ,u_{r}]$ be a ring of polynomials.
A hypermatrix $\CA=(f_{i_{1}, \ldots ,i_{t}})_{1\leq i_{j}\leq n_{j}, \, j=1, \ldots ,t}$, where all $f_{i_{1}, \ldots , i_{t}} \in K[u_{1}, \ldots ,u_{r}]_{1}$, is called a ``weak generic hypermatrix of indeterminates'' (or briefly ``weak generic hypermatrix'') if:
\begin{enumerate}
\item all the entries  of $\CA$ belong to $\{ u_{1}, \ldots ,u_{r}\}$;
\item\label{3} there exists an entry $f_{i_{1}, \ldots , i_{t}}$ such that $f_{i_{1}, \ldots , i_{t}}\neq f_{k_{1}, \ldots , k_{t}}$ for all $(k_{1}, \ldots , k_{t})\neq (i_{1}, \ldots , i_{t})$, $1\leq k_{j}\leq n_{j}, \, j=1, \ldots ,t$;
\item the ideals of $2$-minors of all sections of $\CA$ are prime ideals.
\end{enumerate}
\end{defi}

\begin{lemma}\label{IJ}
Let  $I,J\subset R=K[u_{1}, \ldots , u_{r}]$ be ideals such that $J=(I, u_{1}, \ldots , u_{q})$ with $q<r$. Let $f\in R$ be a polynomial independent of $u_{1}, \ldots , u_{q}$ and such that $I:f=I$. Then $J:f=J$.
\end{lemma}

\begin{proof} We need to prove that if $g\in R$ is such that $fg\in J$, then $g\in J$.

Any polynomial $g\in R$ can be written as $g=g_{1}+g_{2}$ where $g_{1}\in (u_{1}, \ldots , u_{q})$ and $g_{2}$ is independent of $u_{1}, \ldots , u_{q}$. Clearly $g_{1}\in J$.
Now $fg_{2}=fg-fg_{1}\in J$ and $fg_{2}$ is independent of $u_{1}, \ldots , u_{q}$. This implies that $fg_{2}\in I$, then $g_{2}\in I\subset J$ because $I:f=I$ by hypothesis. Therefore $g=g_{1}+g_{2}\in J$.
\end{proof}

Now we can state the main proposition of this section. The proof that we are going to exhibit follows the ideas the proof of Theorem 1.5 in \cite{Ha1}, where the author proves that the ideal generated by 
$2$-minors of a generic hypermatrix of indeterminates is prime. In the same proposition (Proposition 1.12) it is proved that also the ideal generated by $2$-minors of a ``weak generic box'' is prime.  We give here an independent proof for weak generic hypermatrix, since it is a more general result; moreover  we do not follow exactly the same lines as in \cite{Ha1}.

\begin{propos}\label{gen1.12}  Let $R=K[u_{1}, \ldots ,u_{r}]$ be a ring of polynomials and let
 $\CA=(f_{i_{1}, \ldots ,i_{t}})_{1\leq i_{j}\leq n_{j}, \, j=1, \ldots ,t}$ be a weak generic hypermatrix as defined in \ref{weak hypermatrix}. Then the ideal $I_{2}(\CA)$ is a prime ideal in $R$.
\end{propos}

\begin{proof} Since $\CA=(f_{i_{1}, \ldots , i_{t}})_{1\leq i_{j}\leq n_{j}, \; j=1, \ldots , t}$ is a weak generic hypermatrix, there exists an entry $f_{i_{1}, \ldots , i_{t}}$ that verifies the item \ref{3}. in Definition \ref{weak hypermatrix}. It is not restrictive to assume that such $f_{i_{1}, \ldots , i_{t}}$ is $f_{1, \ldots , 1}$.

Let $F,G\in R$ s.t. $FG\in I_{2}(\CA)$. We want to prove that either $F\in I_{2}(\CA)$ or $G\in I_{2}(\CA)$. Let $Z=\{f_{1, \ldots , 1}^{k}\; | \; k\geq 0\}\subset R$ and let $R_{Z}$ be the localization of $R$ at $Z$. Let also $\varphi : R \rightarrow R_{Z}$ such that  
$$\varphi(f_{j_{1}, \ldots , j_{t}})=\frac{f_{j_{1}, 1, \ldots , 1}\cdots f_{1,\ldots , \, 1, j_{t}}}{f_{1, \ldots , 1}^{t-1}},$$
$\varphi(K)=K$ and $\varphi (u_{i})=u_{i}$ for $u_{i}\in \{u_{1}, \ldots , u_{r}\}\backslash \{f_{i_{1}, \ldots , i_{t}}\, | \, 1\leq i_{j}\leq n_{j}, \, j=1, \ldots , t\}$.
Clearly $\varphi(m)=0$ for all $2$-minors $m$ of $\CA$. Hence $\varphi(I_{2}(\CA))=0$.
Since $F(\ldots , f_{j_{1}, \ldots , j_{t}}, \ldots )G(\ldots , f_{j_{1}, \ldots , j_{t}}, \ldots )\in I_{2}(\CA)$ then $F(\ldots , \varphi(f_{j_{1}, \ldots , j_{t}}), \ldots )\cdot G(\ldots , \varphi(f_{j_{1}, \ldots , j_{t}}),\ldots )=0_{R_{Z}}$. The localization $R_{Z}$ is a domain because $R$ is a  domain, thus either $F(\ldots , \varphi(f_{j_{1}, \ldots , j_{t}}), \ldots )=0_{R_{Z}}$, or $G(\ldots , \varphi(f_{j_{1}, \ldots , j_{t}}), \ldots )=0_{R_{Z}}$. Suppose that $F\left(\ldots , \frac{f_{j_{1}, 1, \ldots , 1}\cdots f_{1,\ldots , \, 1, j_{t}}}{f_{1, \ldots , 1}^{t-1}},\ldots\right)=0_{R_{Z}}$. We have
\begin{equation}\label{F+H}F(\ldots , f_{j_{1}, \ldots , f_{j_{t}}}, \ldots )=F\left(\ldots , \frac{f_{j_{1}, 1, \ldots , 1}\cdots f_{1,\ldots , \, 1, j_{t}}}{f_{1, \ldots , 1}^{t-1}},\ldots\right) + H,\end{equation}
where $H$ belongs to the ideal $(f_{j_{1}, \ldots , j_{t}}f_{1, \ldots , 1}^{t-1}-f_{j_{1}, 1 \ldots , 1}\cdots f_{1, \ldots , 1, j_{t}})_{1\leq j_{k}\leq n_{j}, \, k=1, \ldots ,t}\subset R_{Z}$.\\
Now let $H_{t-1}=f_{j_{1}, \ldots , j_{t}}f_{1, \ldots , 1}^{t-1}-f_{j_{1}, 1 \ldots , 1}\cdots f_{1, \ldots , 1, j_{t}}$. Then
$$H_{t-1}=f_{1_{1},j_{2}, \ldots , j_{t}}f_{ j_{1},1,\ldots ,1 }f_{j_{1}, \ldots , j_{t}}^{t-2}+(f_{1, \ldots ,1}f_{j_{1}, \ldots ,j_{t}}-f_{1, j_{2}, \ldots , j_{t}}f_{j_{1}, 1 \ldots , 1})f_{j_{1}, \ldots ,j_{t}}^{t-2}-$$
$$-f_{1, j_{2}, \ldots , j_{t}}f_{j_{1}, 1, j_{3}, \ldots , j_{t}}\cdots f_{j_{1}, \ldots , j_{t-1},1}\equiv_{I_{2}(\CA)}$$
$$f_{1,j_{2},\ldots , j_{t}}f_{j_{1},1, \ldots , 1}f_{1, \ldots , 1}^{t-2}-f_{1, j_{2},\ldots ,j_{t}}f_{j_{1},1,j_{3}, \ldots ,j_{t}}\cdots f_{j_{1}, \ldots , j_{t-1},1}=H_{t-2}.$$
Proceeding analogously for $H_{t-2}, \ldots , H_{1}$, it is easy to verify that $H_{t-1}\in I_{2}(\CA)$. Hence $H$ belongs to the ideal of $R_{Z}$ generated by $I_{2}(\CA)$. This fact, together with (\ref{F+H}), implies that also $F$ belongs to the ideal of $R_{Z}$ generated by $I_{2}(\CA)$. Therefore we obtained that if $\varphi (F)=0_{R_{Z}}$, then there exists $\nu>0$ such that 
\begin{equation}\label{fF}f_{1,\ldots ,1}^{\nu}F(\ldots , f_{j_{1}, \ldots ,j_{t}}, \ldots )\in I_{2}(\CA)\subset R.\end{equation}

Now we want to prove that if there exists $\nu>0$ such that $f_{1,\ldots ,1}^{\nu}F(\ldots , f_{j_{1}, \ldots ,j_{t}}, \ldots )\in I_{2}(\CA)$, then $F\in I_{2}(\CA)$.
Analogously as it is done in the proof of Lemma 1.4 in \cite{Ha1}, we will use a triple induction: first on the dimension $t$ of the hypermatrix $\CA$, then on  $\sum_{j=1}^{t}n_{j}$, and finally on $\deg (F)$.
\begin{description}
\item[Induction on $\mathbf{t}$.] For $t=2$ our goal is proved in Lemma 3 of \cite{Sh}. Assume that $t>2$ and that the induction hypothesis holds for any weak generic hypermatix of size lower than $t$.
\item[Induction on $\mathbf{\sum_{j=1}^{t}n_{j}}$.] If $n_{j}=1$ for at least one $j\in \{1, \ldots , t\}$, then $\CA$ is a hypermatrix of order $(t-1)$, so the result is true for the induction hypothesis on $t$. Assume that $n_{j}\geq 2$ for all $j=1, \ldots , t$ and that the induction hypothesis holds for smaller values of $\sum_{j=1}^{t}n_{j}$.
\item[Induction on $\mathbf{\mathrm{\mathbf{deg}}(F)}$.] If $\deg(F)=0$, since $\varphi (F)=0_{R_{Z}}$, we have $F=0\in I_{2}(\CA)$. Then let $\deg(F)>0$ and assume that the induction hypothesis holds for polynomials of degree lower than $\deg(F)$.
\end{description}
In \cite{Ha1}, Corollary 1.1.1, it is proved that $(I_{2}(\CA), f_{n_{1},\ldots , n_{t}})=\cap_{l=1}^{t}I_{l}$ where $\CA_{l}$ is the hypermatrix $(f_{i_{1}, \ldots , i_{t}})_{i_{l}<n_{l}}$, and $I_{l}:=(I_{2}(\CA_{l}), \{f_{i_{1}, \ldots , i_{t}}\; | \; i_{l}=n_{l}\})$ . Clearly $I_{2}(\CA)\subseteq (I_{2}(\CA), f_{n_{1}, \ldots , n_{t}})$. By (\ref{fF}), we have that $f_{1, \ldots , 1}^{\nu}F\in I_{2}(\CA)$. Hence, by Corollary 1.1.1 in \cite{Ha1}, $f_{1, \ldots , 1}^{\nu}F \in I_{l}$ for all $l=1, \ldots , t$.
We can apply here the induction hypotheses on $t$ and on $\sum_{j=1}^{t}n_{j}$, hence $I_{2}(\CA_{l}): f_{1, \ldots , 1}^{\nu}=I_{2}(\CA_{l})$. Now, by Lemma \ref{IJ}, $I_{l}:f_{1\ldots , , 1}^{\nu}=I_{l}$, i.e. $F\in \cap_{l=1}^{t}I_{l}=(I_{2}(\CA), f_{n_{1}, \ldots , n_{t}})$.
Hence we can write $F=F_{1}+F_{2}$ where $F_{1}\in I_{2}(\CA)$ and $F_{2}\in (f_{n_{1}, \ldots , n_{t}})$, that is to say $F=F_{1}+f_{n_{1}, \ldots , n_{t}}\tilde F_{2}$ with $\deg(\tilde F_{2})<\deg (F)$. Obviously $f_{1, \ldots , 1}^{\nu}f_{n_{1}, \ldots , n_{t}}\tilde F_{2}=f_{1, \ldots , 1}^{\nu}F-f_{1, \ldots , 1}^{\nu}F_{1}\in I_{2}(\CA)$. \\
Let's notice that we checked that, since $\varphi (f_{n_{1}, \ldots , n_{t}})\neq 0_{R_{Z}}$, for any form $K$ for which $f_{n_{1}, \ldots , n_{t}}K\in I_{2}(\CA)$ there exists $\mu>0$ such that $f_{1, \ldots , 1}^{\mu}K\in I_{2}(\CA)$; if we apply this to $K=f_{1, \ldots , 1}^{\nu}\tilde{F}_{2}$, we get that $f_{1, \dots , 1}^{\nu + \mu}\tilde{F}_{2}\in I_{2}(\CA)$ for some $\mu>0$. Now we deduce that there exists $\mu > 0$ s. t. $f_{1, \ldots , 1}^{\nu+\mu}\tilde F_{2}\in I_{2}(\CA)$. Now, by induction hypothesis on the degree of $F$, we have that $\tilde F_{2}\in I_{2}(\CA)$. Therefore $F\in I_{2}(\CA)$.
\end{proof}

\subsection{Ideals of Segre -Veronese varieties}\label{Ideal Segre veronese}

Since a Segre-Veronese variety is given set-theoretically by the $2$-minors of an $(\underline{n}, \underline{d})$-symmetric hypermatrix of indeterminates (see (\ref{SegVerSet})), if we prove that any  $(\underline{n}, \underline{d})$-symmetric hypermatrix of indeterminates is weak generic, we will have, as a consequence of Proposition \ref{gen1.12}, that its $2$-minors are a set of generators for the ideals of  Segre-Veronese varieties.
\\
\\
\Obs
If $\CA=(a_{i_{1}, \ldots , i_{d}})_{1\leq i_{j}\leq n; \, j=1, \ldots , d}$ is a supersimmetric hypermatrix of size $\underbrace{n\times \cdots \times n}_{d}$, then also a $k$-th section $\CA_{i_{k}}^{(l)}$ of $\CA$ is  a supersymmetric hypermatrix of size $\underbrace{n\times \cdots \times n}_{d-1}$.

In fact, since $\CA$ is supersymmetric, then $a_{i_{1}, \ldots , i_{d}}=a_{i_{\sigma(1)}, \ldots , i_{\sigma(d)}}$ for  all $\sigma\in \G_{d}$.
The section $\CA_{i_{k}}^{(l)}$ is obtained from $\CA$ by imposing $i_{k}=l$. Therefore $\CA_{i_{k}}^{(l)}=(a_{i_{1}, \ldots , i_{k}=l , \ldots i_{d}})$ is such that
$a_{i_{1}, \ldots , i_{k}=l , \ldots i_{d}}=a_{i_{\sigma(1)}, \ldots ,i_{\sigma(k)}=l,\ldots , i_{\sigma(d)}}$, for all $\sigma \in \G_{d}$ such that $\sigma(k)=l$, hence such $\sigma$'s can be viewed as elements of  the permutation group of the set $\{1, \ldots , l-1, l+1, \ldots ,d\}$ that is precisely $\G_{d-1}$.
\\
\\
\Obs If $[T]\in Y_{n-1,d}$, then a hypermatrix obtained  as a section of the hypermatrix representing $T$, can be  associated to a tensor $T'$ such that $[T']\in Y_{n-1,d-1}$.

\begin{theorem}\label{theor segre ver} Let $\underline{n}=(n_{1}, \ldots , n_{t})$ and $\underline{d}=(d_{1}, \ldots , d_{t})$. Let  $H_{i}\subset V_{i}^{\otimes d_{i}}$ be the subspace of supersymmetric tensors of $V_{i}^{\otimes d_{i}}$ for $i=1, \ldots , t$ and let $R_{[\underline{n}, \underline{d}]}$ be the ring of coordinates of $\PP {}(H_{1}\otimes \cdots \otimes H_{t})\subset \PP {}(V_{1}^{\otimes d_{1} }\otimes \cdots \otimes V_{t}^{\otimes d_{t}})$ defined in Definition \ref{Sd1dt}.
If $\CA$ is a generic $(\underline{n}, \underline{d})$-symmetric hypermatrix of $R_{[\underline{n}, \underline{d}]}$, then $\CA$ is a weak generic  hypermatrix and the ideal of the Segre-Veronese variety ${\cal S}_{d_{1}, \ldots , d_{t}}(V_{1}\otimes \cdots \otimes V_{t})$ is 
$$I({\cal S}_{d_{1}, \ldots , d_{t}}(V_{1}\otimes \cdots \otimes V_{t}))=I_{2}(\CA)\subset R_{[\underline{n}, \underline{d}]}$$
with $d_{i}>0$ for $i=1, \ldots, t$.
\end{theorem}

\begin{proof} 
The proof is by induction on $\sum_{i=1}^{t}d_{i}$.

The case $\sum_{i=1}^{t}d_{i}=1$ is not very significant because if $\dim(V_{1})=n_{1}$, so ${\cal S}_{1}(V_{1})=Y_{n_{1}-1,1}= \PP {}(V_{1})$, then $I({\cal S}_{1}(V_{1}))=I(\PP {}(V))$ i.e. the zero ideal (in fact the  $2$-minors of $\CA$ do not exist).

If $\sum_{i=1}^{t}d_{i}=2$ the two possible cases for the Segre-Veronese varieties are either ${\cal S}_{2}(V_{1})$ or ${\cal S}_{1,1}(V_{1},V_{2})$. Clearly, if $\dim(V_{1})=n_{1}$, then ${\cal S}_{2}(V_{1})=Y_{n_{1}-1,2}$ is Veronese variety and the theorem holds because of (\ref{veronese}). Analogously ${\cal S}_{1,1}(V_{1},V_{2})=Seg(V_{1}\otimes V_{2})$ and again the theorem is known to be true (\cite{Ha1}).

Assume that the theorem holds for every $(\underline{n},\underline{d})$-symmetric hypermatrix with $\sum_{i=1}^{t}d_{i}\leq r-1$. Then, by Proposition \ref{gen1.12}, the ideal generated by the $2$-minors of such an $(\underline{n}, \underline{d})$-symmetric hypermatrix is a prime ideal.

Now, let $\CA$ be an $(\underline{n},\underline{d})$-symmetric hypermatrix with $\sum_{i=1}^{t}d_{i}=r$.  The first two properties that characterize a weak generic hypermatrix (see Definition \ref{weak hypermatrix}) are immediately verified for $\CA$.  For the third one we have to check that the ideals of the $2$-minors of all sections $\CA_{i_{p,q}}^{(l)}$ of $\CA$ are prime ideals. \\
If we prove that $\CA_{i_{p,q}}^{(l)}$ represents an $(\underline{n}, \underline{d}')$-symmetric hypermatrix (with $\underline{d}'=(d_{1}, \ldots , d_{p}-1 , \ldots , d_{t}))$) we will have, by induction hypothesis, that $\CA_{i_{p,q}}^{(l)}$ is a weak generic hypermatrix and hence its $2$-minors generate a prime ideal.
\\
The hypermatrix  $\CA=(a_{i_{1,1}, \ldots , i_{1,d_{1}}; \ldots ; i_{t,1} ,\ldots i_{t,d_{t}}})_{1\leq i_{j,k}\leq n_{j}, \, k=1, \ldots , d_{j}, \, j=1, \ldots , t}$ is $(\underline{n},\underline{d})$-symmetric, hence, by definition, $a_{i_{1,1}, \ldots , i_{1,d_{1}}; \ldots ; i_{t,1}, \ldots , i_{t,d_{t}}}=a_{i_{\sigma_{1}(1,1)}, \ldots , i_{\sigma_{1}(1,d_{1})}; \ldots ; i_{\sigma_{t}(t,1)}, \ldots , i_{\sigma_{t}(t,d_{t})}}$ for all permutations $\sigma_{j}\in \G(j,d_{j})$ where $\G(j,d_{j})$ is the permutation group on $\{(j,1), \ldots , (j,d_{j})\}$ for all $j=1, \ldots , t$. 
\\
The hypermatrix $\CA_{i_{p,q}}^{(l)}=(a_{i_{1,1}, \ldots , i_{1,d_{1}}; \ldots  ,i_{p,q}=l,\ldots ;i_{t,1}, \ldots , i_{t,d_{t}}})$, obtained from $\CA$ by imposing $i_{p,q}=l$, is $(\underline{n}, \underline{d}')$-symmetric because
$$a_{i_{1,1}, \ldots , i_{1,d_{1}}; \ldots , i_{p,q}=l,\ldots ;  i_{t,1}, \ldots , i_{t,d_{t}}}=a_{i_{\sigma_{1}(1,1)}, \ldots , i_{\sigma_{1}(1,d_{1})}; \ldots ,i_{\sigma_{p}(p,1)}, \ldots , i_{p,q}=l, \ldots  i_{\sigma_{p}(p,d_{p})}; \ldots ; i_{\sigma_{t}(t,1)}, \ldots , i_{\sigma_{t}(t,d_{t})}}$$ 
for all $\sigma_{j}\in \G(j,d_{j})$, $j=1, \ldots , \hat p , \ldots , t$, and for $\sigma_{p}\in \G(p,d_{p}-1)$, where $\G(p,d_{p}-1)$ is the permutation group  on the set of indices $\{(p,1), \ldots , \widehat{(p,q)}, \ldots , (p,d_{p})\}$ (this is a consequence of the first Remark of this section). Hence $I_{2}({\CA}_{i_{p,q}}^{(l)})$ is prime by induction, and $\CA$ is weak generic, so also $I_{2}(\CA)$ is prime.

Since  by definition ${\cal S}_{d_{1}, \ldots , d_{t}}(V_{1}\otimes \cdots \otimes V_{t})=\PP {}(H_{1}\otimes \cdots \otimes H_{t})\cap Seg(V_{1}\otimes \cdots \otimes V_{t})$, we have that $I_{2}(\CA)$ is a set of equations for ${\cal S}_{d_{1}, \ldots ,d_{t}}(V_{1}\otimes \cdots \otimes V_{t})$ (see (\ref{SegVerSet})), hence, because of the primeness of $I_{2}(\CA)$ that we have just proved, $I_{2}(\CA)\subset R_{[\underline{n}, \underline{d}]}$ is the ideal of ${\cal S}_{d_{1}, \ldots ,d_{t}}(V_{1}\otimes \cdots \otimes V_{t})$.
\end{proof}

\section{Projections of Veronese surfaces}\label{PoVs}

In this section we want to use the tool of weak generic hypermatrices in order to prove that the ideal of a projection of a Veronese surface $Y_{2,d}\subset \PP {{d+2\choose d}-1}$ from a finite number $s \leq {d\choose 2}$ of general points on it is the prime ideal defined by the order $2$-minors of some particular tensor.

In \cite{Ha1}  the case in which $s$ is a binomial number (i.e. $s={t+1\choose 2}$ for some positive integer $t\leq d-1$) is done.

In this section we try to extend that result to a projection of a Veronese surface from any number $s\leq {d\choose 2}$ of general points.

Notice that in \cite{Gi} and in  \cite{GL} the authors study the projection of Veronese surfaces $Y_{2,d}$ from $s={d\choose 2}+k$ general points, $0\leq k \leq d$, for some non negative integer $k$, (this corresponds to the case of a number of points between the two consecutive binomial numbers ${d\choose 2}$ and ${d+1\choose 2}$).
\\
\\
Let $Z=\{P_{1}, \ldots , P_{s}\}\subset \PP 2$ be a set of general points in $\PP 2$, where $s={t+1\choose 2}+k\leq {d\choose 2}$ with $0<t\leq d-1$ and $0\leq k \leq t$ (actually we may assume $t\leq d-2$ because the case $t=d-1$ and $k=0$ corresponds to the known case of the ``Room Surfaces'' - see \cite{GG}). Let $J\subset S= K[w_{1},w_{2},w_{3}]$ be the ideal $J=I(Z)$, i.e. $J=\wp_{1} \cap \cdots \cap \wp_{s}$ with $\wp_{i}=I(P_{i})\subset S$ prime ideals for $i=1, \ldots , s$. 

Let  $J_{d}$ be the degree $d$ part of the ideal $J$  and let $Bl_{Z}(\PP 2)$ be the blow up of $\PP 2$ at  $Z$. Since $d\geq t+1$,  the linear system of the strict transforms of the curves defined by $J_{d}$, that we indicate with $\tilde{J}_{d}$,  is very ample.
If $\varphi_{J_{d}}: \PP 2 \dashrightarrow \PP {{d+2\choose 2}-s-1}$ is the rational morphism associated to $J_{d}$ and if $\varphi_{\tilde{J}_{d}}:Bl_{Z}(\PP 2) \rightarrow \PP {{d+2\choose 2}-s-1}$ is the morphism associated to $\tilde{J}_{d}$, the variety $X_{Z,d}$ we want to study is $\overline{Im(\varphi_{J_{d}})}=Im(\varphi_{\tilde{J}_{d}})$.
 This variety can also be viewed as the projection of the Veronese surface $Y_{2,d}\subset \PP {{d+2\choose 2}-1}$ from $s$ general points on it.
 
The first thing to do is to describe $J_{d}$ as vector space.

\subsection{The ideal of general points in the projective plane}\label{4.1}

There is a classical result, Hilbert-Burch Theorem (see, for instance, \cite{CGO}), that gives a description of the generators of $J$. I.e. the ideal $J\subset S=K[w_{1},w_{2},w_{3}]$ is generated by $t-k+1$ forms $F_{1}, \ldots , F_{t-k+1}\in S_{t}$ and by $h$ forms $G_{1}, \ldots , G_{h}\in S_{t+1}$ where $h=0$ if $0\leq k < t/2$ and $h=2k-d$ if $t/2\leq k\leq t$. What follows now is the constructions of the $F_{j}$'s and the $G_{i}$'s (the same description is presented in \cite{GL}).

\begin{description}
\item[If $\mathbf{t/2\leq k\leq t}$,] for a general choice of points $P_{1}, \ldots , P_{s}$, the generators of $J$ can be chosen to be the maximal minors of:
\begin{equation}\label{Ltpiccolo}{\mathcal{L}}:= \left( \begin{array}{ccccccc}
L_{1,1}& \cdots &L_{1,2k-t}&Q_{11}&\cdots &Q_{1,t-k+1}\\
\vdots&&\vdots&\vdots&&\vdots \\
L_{k,1}&\cdots &L_{k,2k-t}&Q_{k,1}&\cdots &Q_{k,t-k+1}
\end{array}\right) \in M_{k,k+1}(S)\end{equation}
where $L_{i,j}\in S_{1}$ and $Q_{h,l}\in S_{2}$ for all $i,h=1, \ldots , k$, $j=1, \ldots , 2k-t$ and $l=1, \ldots , t-k+1$.\\
The forms $F_{j}\in S_{t}$ are the minors of $\mathcal{L}$ obtained by deleting the $2k-t+j$-th column, for $j=1, \ldots , t-k+1$; the forms $G_{i}\in S_{t+1}$ are the minors of $\mathcal{L}$ obtained by deleting the $i$-th column, for $i=1, \ldots , 2k-t$.\\
The degree $(t+1)$ part of the ideal $J$ is clearly $J_{t+1}=<w_{1}F_{1}, \ldots , w_{3}F_{t-k+1}, G_{1}, \ldots , G_{2k-t}>$. If we set $\tilde{G}_{i,j}=w_{i}F_{j}$ for $i=1,2,3$, $j=1, \ldots , t-k+1$ we can write:
$$J_{t+1}=<\tilde{G}_{1,1}, \ldots , \tilde{G}_{3,t-k+1}, G_{1}, \ldots , G_{2k-t}>.$$
Notice that  $w_{1}F_{1}=\tilde{G}_{1,1}, \ldots , w_{3}F_{t-k+1}=\tilde{G}_{3,t-k+1}$ are linearly independent (see, for example, \cite{CGO}). 

\item[If $\mathbf{0\leq k <t/2}$,] then $J$ is generated by maximal minors of:
\begin{equation}\label{Ltgrande}{\mathcal{L}}:=\left( \begin{array}{ccccc}
Q_{1,1} &\cdots&\cdots&\cdots &Q_{1,t-k+1}\\
\vdots &&&& \vdots \\
Q_{k,1}&\cdots &\cdots&\cdots& Q_{k,t-k+1}\\
L_{11}&\cdots &\cdots&\cdots& L_{1,t-k+1}\\
\vdots &&&&\vdots \\
L_{t-2k, 1} &\cdots&\cdots&\cdots&L_{t-2k,t-k+1}
\end{array}\right)\in M_{t-k,t-k+1}(S)\end{equation}
where $L_{i,j}\in S_{1}$ and $Q_{h,l}\in S_{2}$ for all $i=1, \ldots , t-2k$, $j,l=1, \ldots , t-k+1$ and $h=1,\ldots ,k$.
\\
The forms $F_{j}\in S_{t}$ are the minors of $\mathcal{L}$ obtained by deleting the $j$-th column for $j=1, \ldots , t-k+1$.\\
Again $J_{t+1}=<w_{1}F_{1}, \ldots , w_{3}F_{t-k+1}>$ but now those generators are not necessarily linearly independent. Using the same notation of the previous case one can write:
$$J_{t+1}=<\tilde{G}_{1,1}, \ldots , \tilde{G}_{3,t-k+1}>.$$
\end{description}

Clearly if $t/2 \leq k \leq t$ then:
\begin{equation}\label{Jd}J_{d}=<\w^{d-t-1}\tilde{G}_{i,j}, \w^{d-t-1}G_{l}>\end{equation}
for $i=1,2,3$, $j=1, \ldots t-k+1$, $l=1, \ldots , 2k-t$ and $\w^{d-t-1}G=\{w_{1}^{d-t-1}G, w_{1}^{t-d-2}w_{2}G, \ldots , w_{3}^{d-t-1}G\}$.
\\
If $0\leq k<t/2$ then:
\begin{equation}\label{Jd2}J_{d}=<\w^{d-t-1}\tilde{G}_{i,j}>\end{equation}
for $i=1,2,3$ and $j=1, \ldots , t-k+1$.

Denote
$$\left\{
\begin{array}{l}
z_{1}:=w_{1}^{d-t-1}, \\ z_{2}:=w_{1}^{t-d-2}w_{2}, \\ \vdots \\ z_{u}:=w_{3}^{t-d-1}
\end{array}
\right.$$
where $u={d-t+1\choose 2}$; or $z_{\underline{\alpha}}$ for $\underline{w}^{\underline{\alpha}}=w_{1}^{\alpha_{1}}w_{2}^{\alpha_{2}}w_{3}^{\alpha_{3}}$, if $\underline{\alpha}=(\alpha_{1},\alpha_{2},\alpha_{3})\in \NN^{3}$,  $|\underline{\alpha}|=d-t-1$ and we assume that the $\underline{\alpha}$'s are ordered by the lexicographic order.

Let $N$ be the number of generators of $J_{d}$, and let $K[\tilde{x}_{h;i,j}, x_{h,l}]$ be a ring of coordinates on $\PP {N-1}$ with $l=1, \ldots , 2k-t$ only if $t/2\leq k \leq t$ (in the other case the variables $x_{h,l}$ do not exist at all) and $h=1, \ldots , u$; $i=1,2,3$; $j=1, \ldots , t-k+1$ in any case. The morphism $\varphi : \PP 2 \setminus Z \rightarrow \PP {N-1}$ such that
$$\varphi ([w_{1},w_{2},w_{3}])=[z_{1}\tilde{G}_{1,1}, \ldots , z_{u}\tilde{G}_{3,t-k+1},z_{1}G_{1}, \ldots ,z_{u}G_{2k-t}] \hbox{, if } t/2\leq k \leq t,$$
or
$$\varphi([w_{1},w_{2},w_{3}])=[z_{1}\tilde{G}_{1,1}, \ldots , z_{u}\tilde{G}_{3,t-k+1}] \hbox{, if } 0\leq k <t/2,$$
gives a parameterization of $X_{Z,d}$ into $\PP {N-1}$. Observe that $X_{Z,d}=\overline{\varphi_{J_{d}}(\PP 2 \setminus Z)}$ is naturally embedded into $\PP {{d+2\choose 2}-s-1}$, because $\dim_{K}(J_{d})={d+2\choose 2}-s$. In terms of the $\tilde{x}_{h;i,j}$'s and the $x_{h,l}$'s, since the parameterization of $X_{Z,d}$ is:
\begin{equation}\label{x}\left\{
\begin{array}{l}
\tilde{x}_{h;i,j}=z_{h}\tilde{G}_{i,j},\\
x_{h,l}=z_{h}G_{l},
\end{array}
\right. \end{equation}
the independent linear relations between the generators of $J_{d}$ will give the subspace $\PP {}(<Im(\varphi_{\tilde{J}_{d}})>)=\PP {{d+2\choose 2}-s-1}$ of $\PP {N-1}$.
The number of such relations has to be $N-{d+2\choose 2}+s$.

If $t/2\leq k\leq t$, the number of generators of $J_{d}$ given by (\ref{Jd}) is ${d-t+2\choose 2}(t-k+1)+{d-t-1+2\choose 2}(2k-t)$; hence there must be ${d-t\choose 2}k$ independent relations between those generators of $J_{d}$.

If $0\leq k<t/2$, the number of generators of $J_{d}$ in (\ref{Jd2}) is ${d-t+2\choose 2}(t-k+1)$, hence there must be ${d-t+1\choose 2}(t-k)-k(d-t)$ independent relations between those generators of $J_{d}$.

There is a very intuitive way of finding exactly those numbers of relations between the generators of $J_{d}$ and this is what we are going to describe (then we will prove that such relations are also independent).

\begin{description}
\item[If $\mathbf{t/2\leq k\leq t}$,] assume that $\underline{\beta}=(\beta_{1}, \beta_{2}, \beta_{3})$ with $|\underline{\beta}|=d-t-2$. The determinant obtained by adding to the matrix $\mathcal{L}$ defined in (\ref{Ltpiccolo}) a row $\left( \begin{array}{cccccc}
\w^{\underline{\beta}}L_{i,1}&\cdots &\w^{\underline{\beta}}L_{i,2k-t}&\w^{\underline{\beta}}Q_{i,1}& \cdots & \w^{\underline{\beta}}Q_{i,t-k+1}
\end{array}\right)$
clearly vanish for all $i=1, \ldots ,k$:
$$\det \left( \begin{array}{cccccc}
\w^{\underline{\beta}}L_{i,1}&\cdots &\w^{\underline{\beta}}L_{i,2k-t}&\w^{\underline{\beta}}Q_{i,1}& \cdots & \w^{\underline{\beta}}Q_{i,t-k+1}\\
&&&{\mathcal{L}}&&
\end{array}\right) =0.$$
Computing those determinants, for $i=1, \ldots ,k$, one gets:
\begin{equation}\label{primopasso}\sum_{r=1}^{2k-t}\w^{\underline{\beta}}L_{i,r}G_{r}+\sum_{p=1}^{t-k+1}\w^{\underline{\beta}}Q_{i,p}F_{p}=0\end{equation}
where the $G_{r}$'s and the $F_{p}$'s are defined as minors of (\ref{Ltpiccolo}).
\\
Since $L_{i,r}\in S_{1}$, there exist some $\lambda_{i,r,l}\in K$, for $i=1, \ldots , k$, $r=1, \ldots , 2k-t$ and $l=1,2,3$, such that 
$$L_{i,r}=\sum_{l=1}^{3}\lambda_{i,r,l}w_{l};$$ 
analogously, since $Q_{i,p}\in S_{2}$, there exist some $\gamma_{i,p,l,h}\in K$, for $i=1, \ldots , k$, $p=1, \ldots , t-k+1$ and $l,h=1,2,3$, such that $$Q_{i,p}=\sum_{l,h=1}^{3}\gamma_{i,p,l,h}w_{l}w_{h}.$$ 
Before rewriting the equations (\ref{primopasso}), observe that
$$Q_{i,p}F_{p}=\left(\sum_{l,h=1}^{3}\gamma_{i,p,l,h}w_{l}w_{h}\right)F_{p}=\sum_{l,h=1}^{3}\gamma_{i,p,l,h}w_{l}\tilde{G}_{h,p},$$ 
and set: 
\begin{itemize}
\item $\mu_{i,\underline{\alpha} , r}=\left\{ \begin{array}{ll}
\lambda_{i,r,l}, & \hbox{if } \underline{\alpha}=\underline{\beta} +\underline{e}_{l},\\
0 & \hbox{otherwise,}
\end{array}\right.$
\\
for $i=1, \ldots , k$; $|\underline{\alpha} |=t-d-1$ and $l=1,2,3$ and where $\underline{e}_{1}=(1,0,0)$, $\underline{e}_{2}=(0,1,0)$ and $\underline{e}_{3}=(0,0,1)$;
\item $\tilde{\mu}_{i,\underline{\alpha} , p,h}=\left\{ \begin{array}{ll}
\gamma_{i,p,l,h}, & \hbox{if } \underline{\alpha}=\underline{\beta} +\underline{e}_{l},\\
0 & \hbox{otherwise,}
\end{array}\right.$
\\
for $i=1, \ldots , k$; $p=1, \ldots , t-k+1$; $l,h=1,2,3$ and $|\underline{\alpha}|=d-t-2$.
\end{itemize}
Therefore the equations (\ref{primopasso}), for $i=1, \ldots , k$, can be rewritten as follows:
\begin{equation}\label{reltraG}
\sum_{\footnotesize{\begin{array}{c}
|\underline{\alpha}|=d-t-1
\\
1\leq r \leq 2k-t
\end{array}}} 
\mu_{i,\underline{\alpha} , r}\w^{\underline{\alpha}}G_{r} +
\sum_{\footnotesize{\begin{array}{c}
|\underline{\alpha}|=d-t-1\\
1\leq p \leq t-k+1 \\
h=1,2,3
\end{array}}}
\tilde{\mu}_{i,\underline{\alpha}, p , h} \w^{\underline{\alpha}}\tilde{G}_{h,p} 
=0,\end{equation}
which, for $i=1, \ldots , k$, in terms of $x_{\underline{\alpha},r}$ and $\tilde{x}_{\underline{\alpha},h,p}$  defined in (\ref{x}) becomes:
\begin{equation}\label{equinx}\tag{$E_{1}$}
\sum_{\footnotesize{\begin{array}{c}
|\underline{\alpha}|=d-t-1
\\
1\leq r \leq 2k-t
\end{array}}} 
\mu_{i,\underline{\alpha} , r}x_{\underline{\alpha}, r} +
\sum_{\footnotesize{\begin{array}{c}
|\underline{\alpha}|=d-t-1\\
1\leq p \leq t-k+1 \\
h=1,2,3
\end{array}}}
\tilde{\mu}_{i,\underline{\alpha}, p , h} \tilde{x}_{\underline{\alpha},h,p} 
=0.
\end{equation}
There are exactly $k$ of such relations for each $\underline{\beta}$ and the number of $\underline{\beta}$'s is ${d-t\choose 2}$. Hence in (\ref{reltraG}) we have found precisely the number of relations between the generators of $J_{d}$ that we were looking for; we need to prove that they are independent.

\item[If $\mathbf{0\leq k<t/2}$,] the way of finding the relations between the generators of $J_{d}$ is completely analogous to the previous one. The only difference is that in this case they come from the vanishing of two different kinds of determinants:
\begin{equation}\label{detL}
\det \left( \begin{array}{ccc}
\w^{\underline{\beta}}L_{i,1} &\cdots &\w^{\underline{\beta}}L_{i,t-k+1}\\
&{\mathcal{L}}&
\end{array}\right) =0
\end{equation}
for $i=1, \ldots  , t-2k$, $|\underline{\beta}|=d-t-1$ and $\mathcal{L}$ defined as in (\ref{Ltgrande});
and
\begin{equation}\label{detQ}
\det \left( \begin{array}{ccc}
\w^{\underline{\beta}'}Q_{j,1} &\cdots &\w^{\underline{\beta}'}Q_{j,t-k+1}\\
&{\mathcal{L}}&
\end{array}\right) =0
\end{equation}
for $j=1, \ldots , k$, $|\underline{\beta}'|=d-t-2$ and $\mathcal{L}$ defined as in (\ref{Ltgrande}).

Proceeding as in the previous case one finds that the relations coming from (\ref{detL}) are of the form
\begin{equation}\label{fromL}\tag{$E$}
\sum_{\footnotesize{\begin{array}{c}
|\underline{\alpha}|=d-t-1
\\
1\leq r \leq t-k-1\\
l,h=1,2,3
\end{array}}}
\tilde{\lambda}_{i,\underline{\alpha} ,r,l}z_{\underline{\alpha}}\tilde{G}_{h, r}=0\end{equation}
for some $\tilde{\lambda}_{i,\underline{\alpha} ,r,l}\in K$ and the number of them is ${d-t+1\choose 2}(t-2k)$. 

The relations coming from (\ref{detQ}) are of the form
\begin{equation}\label{fromQ}\tag{$EE$}
\sum_{\footnotesize{\begin{array}{c}
|\underline{\alpha}|=d-t-1\\
1\leq r \leq t-k+1 \\
l,h=1,2,3
\end{array}}}
\tilde{\mu}_{i,\underline{\alpha}, r, l}z_{\underline{\alpha}}\tilde{G}_{h, r}=0\end{equation}
for some $\tilde{\mu}_{i,\underline{\alpha}, r, l}\in K$ and the number of them is ${d-t\choose 2}k$.

The equations (\ref{fromL}) and (\ref{fromQ}) allow to observe that $X_{Z,d}$ is contained in the projective subspace of $\PP {N-1}$ defined by the following linear equations in the variables $\tilde{x}_{\underline{\alpha} , h,r}$:

\begin{equation}\label{E2}\tag{$E_{2}$}
\left\{ \begin{array}{l}
\sum_{\footnotesize{\begin{array}{c}
|\underline{\alpha}|=d-t-1
\\
1\leq r \leq t-k-1\\
l,h=1,2,3
\end{array}}}
\tilde{\lambda}_{i,\underline{\alpha} ,r,l}\tilde{x}_{\underline{\alpha} ;h, r}=0\\
\sum_{\footnotesize{\begin{array}{c}
|\underline{\alpha}|=d-t-1\\
1\leq r \leq t-k+1 \\
l,h=1,2,3
\end{array}}}
\tilde{\mu}_{i,\underline{\alpha}, r, l}\tilde{x}_{\underline{\alpha} ;h, r}=0
\end{array}
\right.
\end{equation}

The number of relations (\ref{E2}) is ${d-t+1\choose 2}(t-2k)+{d-t\choose 2}k$, that is exactly the number of independent relations we expect in the case $0\leq k<t/2$.
\end{description}

Now we have to prove that the relations (\ref{equinx}), respectively (\ref{E2}),  are independent.
\\
\\
\Not
Let $M$ be the matrix of order $\left( {d-t\choose 2}k\right) \times \left({d-t+1\choose 2}(2t-k+3)\right)$ given by the $\mu_{i,\underline{\alpha} , r}$ and the $\tilde{\mu}_{i,\underline{\alpha} , p,h}$ appearing in all the equations (\ref{equinx}). We have already observed that there exists an equation of type (\ref{equinx}) for each  multi-index over three variables  $\underline{\beta}$ of weight $|\underline{\beta}|=d-t-2$, and for each $i=1,\ldots , k$. We construct the matrix $M$ by blocks $M_{\underline{\beta}, \underline{\alpha}}$ (the triple multi-index $\underline{\alpha}$ is such that $|\underline{\alpha}|=d-t-1$): 
\begin{equation}\label{M}M=(M_{\underline{\beta},\underline{\alpha}})_{|\underline{\beta}|=d-t-2, |\underline{\alpha}|=d-t-1}\end{equation}
and the orders on the $\underline{\beta}$'s and the $\underline{\alpha}$'s are the respective decreasing lexicographic orders.
For each fixed $\underline{\beta}$ and $\underline{\alpha}$, the block $M_{\underline{\beta}, \underline{\alpha}}$ is the following matrix:
$$M_{\underline{\beta},\underline{\alpha}}=\left(\begin{array}{cccccc}
\mu_{1,\underline{\alpha}, 1} & \cdots &\mu_{1,\underline{\alpha} , 2k-t}& \tilde{\mu}_{1,\underline{\alpha} , 1,1} & \cdots  & \tilde{\mu}_{1,\underline{\alpha} ,t-k+1,3}  \\
\vdots &&\vdots &\vdots &&\vdots \\
\mu_{k,\underline{\alpha}, 1} & \cdots &\mu_{k,\underline{\alpha} , 2k-t}& \tilde{\mu}_{k,\underline{\alpha} , 1,1} & \cdots  & \tilde{\mu}_{k,\underline{\alpha} , t-k+1,3}
\end{array}\right).$$

Analogously we construct the matrix $N$ of order $\left({d-t+1\choose 2}(t-2k)+{d-t\choose 2}k \right)\times\left( 3{t-d+1\choose 2}(t-k+1)\right)$:
\begin{equation}\label{N}N:=\left(\begin{array}{c}
N_{\underline{\beta}, \underline{\alpha}}\\
N_{\underline{\beta}', \underline{\alpha}}
\end{array}\right)_{|\underline{\alpha}|=|\underline{\beta}|=d-t-1,\, |\underline{\beta}' |=d-t-2}
\end{equation}
where
$$N_{\underline{\beta}, \underline{\alpha}}:=\left( \begin{array}{ccc}
\tilde{\lambda}_{1,\underline{\alpha} , 1,1} &\cdots & \tilde{\lambda}_{1,\underline{\alpha} , t-k-1,3}\\
\vdots && \vdots \\
\tilde{\lambda}_{t-2k,\underline{\alpha} , 1,1}&\cdots &\tilde{\lambda}_{t-2k,\underline{\alpha} , t-k-1,3}
\end{array}\right) \; \; \hbox{ and }\; \; N_{\underline{\beta}',\underline{\alpha}}:=\left( \begin{array}{ccc}
\tilde{\mu}_{1,\underline{\alpha} 1, 1} &\cdots & \tilde{\mu}_{1,\underline{\alpha} t-k+1, 3}\\
\vdots && \vdots \\
\tilde{\mu}_{k,\underline{\alpha} 1, 1}&\cdots &\tilde{\mu}_{k,\underline{\alpha} t-k+1, 3}
\end{array}\right)$$
where the $\tilde{\lambda}_{i,\underline{\alpha} , r,l}$'s and the $\tilde{\mu}_{i , \underline{\alpha} , r, l}$'s are those appearing in 
(\ref{fromL}) and in (\ref{fromQ}) respectively.

\begin{propos} 
The matrices $M$ and $N$ defined in (\ref{M}) and (\ref{N}), respectively, are of maximal rank.
\end{propos}

\begin{proof} Without loss of generality we may assume that  $P=[0,0,1]\notin Z$ and that $F_{1}$ (i.e. the first minor of the matrix $\mathcal{L}$ defined either in (\ref{Ltpiccolo}) or in (\ref{Ltgrande})) does not vanish at $P$. 

 For the $M$ case, one can observe that every time $\underline{\alpha}\neq \underline{\beta} + \underline{e}_{l}$, $l=1,2,3$, the block $M_{\underline{\beta},\underline{\alpha}}$ is identically zero, and we denote $M_{\underline{\beta} , \underline{\beta} +\underline{e}_{l}}$ with $A_{l}$ for $l=1,2,3$.

Consider $\tilde{M}$ the maximal square submatrix of $M$ obtained by deleting the last columns of $M$ (recall that we have ordered both the columns and the rows of $M$ with the respective decreasing lexicographic orders). 

All the blocks $M_{\underline{\beta} ,\underline{\alpha}}$ on the diagonal of $\tilde M$ are such that the position of $\underline{\beta}$ is the same position of $\underline{\alpha}$ in their respective decreasing lexicographic orders. Since $|\underline{\beta}|=|\underline{\alpha}|-1$, then the  blocks appearing on the diagonal of $\tilde M$ are $M_{\underline{\beta} , \underline{\beta} + \underline{e}_{1}}=A_{1}$ for all $\underline{\beta}$'s.
\\
If $\underline{\beta}=(\beta_{1},\beta_{2},\beta_{3})$ and $\underline{\alpha}=(\alpha_{1},\alpha_{2}, \alpha_{3})$, the blocks $M_{\underline{\beta} , \underline{\alpha}}$ under the diagonal are all such that $\beta_{1}<\alpha_{1}-2$, hence they are all equal to zero.
\\
This is clearly sufficient to prove that $\tilde M$ has maximal rank; then $M$ has maximal rank too.
 
The $N$ case is completely analogous.
\end{proof}

With this discussion we have proved the following:

\begin{propos} The coordinates of the points in $X_{Z,d}\subset \PP {N-1}=\PP {}((K[\tilde{x}_{h;i,j},x_{h,l}]_{1})^{*})$ satisfy either the equations (\ref{equinx}) if $t/2\leq k \leq t$, or (\ref{E2}) if $0\leq k <t/2$. Moreover the relations (\ref{equinx}), respectively (\ref{E2}), are linearly independent.
\end{propos}

\Obs There exist other linear relations between the $\tilde{x}_{\underline{\alpha};i,j}$'s and the $x_{\underline{\alpha},l}$ coming from the fact that $w_{i}\tilde{G}_{h,j}=w_{h}\tilde{G}_{i,j}$ for $i,h=1,2,3$ and all $j$'s. If  we denote $z_{\underline{\beta}+\underline{e}_{i}}=\w^{\underline{\beta}}w_{i}$ (with $|\underline{\beta}|=d-t-2$), we have that $z_{\underline{\beta} + \underline{e}_{i}}\tilde{G}_{h,j}=z_{\underline{\beta}+\underline{e}_{h}}\tilde{G}_{i,j}$,
 that is equivalent to: 
$$
\tilde{x}_{\underline{\beta}+\underline{e}_{i}; h,j}=\tilde{x}_{\underline{\beta} +\underline{e}_{h}; i,j}.$$
 The proposition just proved and the fact that the span $<Im(\varphi_{\tilde{J}_{d}})>$ has the same dimension of the subspaces of $\PP N$ defined by  either (\ref{equinx}) or by (\ref{E2}), imply that those relations are linear combinations of either the  (\ref{equinx}), or the (\ref{E2}).
\\
\\
Now the study moves from the linear dependence among generators of $J_{d}$ to the dependence in higher degrees.

\subsection{Quadratic relations}

\Obs 
\begin{enumerate}
\item\label{prova} Let  $X:=(\tilde{x}_{h;i,j},x_{h,l})_{h; i,j,l}$ be the matrix whose entries are the variables of the coordinate ring $K[\tilde{x}_{h;i,j},x_{h,l}]_{1}$
where the index $h=1, \ldots , {d-t+1\choose 2}$ indicates the rows of $X$, and the indicies $(i,j,l)$ indicate the columns and are ordered via the lexicographic order, $i=1,2,3$, $j=1\ldots , t-k+1$, $l=1, \ldots ,2k-t$ (when it occurs). 
\\
The  $2$-minors of $X$  are annihilated by points of $X_{Z,d}$. Denote this set of equations with (XM).
\item The $z_{i}$'s satisfy the equations of the Veronese surface $Y_{2,d-t-1}$, i.e. the $2$-minors of the following catalecticant matrix:
\begin{equation}\label{C}C:=\left(\begin{array}{ccccc}
z_{1} & z_{2} &z_{3}& \cdots & z_{u-2}\\
z_{2}&z_{4}&z_{5}&\cdots&z_{u-1}\\
z_{3}&z_{5}&z_{6}&\cdots &z_{u}
\end{array}\right) \end{equation}
with $u={d-t+1\choose 2}$.
\\
Multiplying $C$ either by $\tilde G_{i,j}$, or by $G_{l}$, for each $i=1,2,3$; $j=1, \ldots , t-k+1$ and $l=1, \ldots , 2k-t$, one obtains either
$$\left( \begin{array}{ccc}
\tilde{x}_{1;i,j} &\cdots&\tilde{x}_{u-2;i,j}\\
\tilde{x}_{2;i,j} &\cdots&\tilde{x}_{u-1;i,j}\\
\tilde{x}_{3;i,j} &\cdots&\tilde{x}_{u;i,j}
\end{array}\right) 
\hbox{, or }
\left(\begin{array}{ccc}
x_{1,l}&\cdots&x_{u-2,l}\\
x_{2,l}&\cdots&x_{u-1,l}\\
x_{3,l}&\cdots&x_{u,l}
\end{array}\right).$$
Therefore on $X_{Z,d}\subset \PP {N-1}$, the coordinates 
$\tilde{x}_{1;i,j}, \ldots , \tilde{x}_{u;i,j}$, for all $i=1,2,3$ and $j=1, \ldots , t-k+1$, or $x_{1,l}, \ldots , x_{u,l}$, for all $l=1, \ldots , 2k-t$, annihilate the $2$-minors of those catalecticant  matrices, respectively. Denote the set of all these equations with (Cat).
\item 
For all $h=1, \ldots , {d-t+1\choose 2} $, on $X_{Z,d}$ we have that $\tilde{G}_{i,j}=\tilde{x}_{h,i,j}/z_{h}$  and  $G_{l}=x_{h,l}/z_{h}$ therefore on $X_{Z,d}\times Y_{2,d-t-1}$ the following system of equations is satisfied for all $h$'s:
\begin{equation}\label{Sh}\tag{$S_{h}$}\left\{ \begin{array}{l}
\tilde{x}_{h;i,j}z_{1}=\tilde{x}_{1;i,j}z_{h}\\
\vdots \\
\tilde{x}_{h;i,j}z_{u}=\tilde{x}_{u;i,j}z_{h}\\
x_{h,l}z_{1}=x_{1,l}z_{h}\\
\vdots \\
x_{h,l}z_{u}=x_{u,l}z_{h}\\
\end{array} \right.
 .\end{equation}
\end{enumerate}

\begin{propos}\label{PQ}
Let $Q:[\tilde{x}_{h;i,j},x_{h,l}]$, $h=1, \ldots , {d-t+1\choose 2}$,  $i=1,2,3$, $j=1\ldots , t-k+1$ and $l=1, \ldots ,2k-t$, such that the equations (XM) are zero if evaluated in $Q$. Then there exists a point $P:[z_{1}, \ldots , z_{u}]\in \PP {u-1}$ such that $P$ and $Q$ satisfy the equations (\ref{Sh}) for all $h$'s.
\end{propos}

\begin{proof}
Since $Q:[\tilde{x}_{1;1,1}, \ldots , x_{{d-t+1\choose 2},2k-t}]$ annihilates all the equations (XM),  the rank of $X$ at $Q$ is $1$, i.e., if we assume that the first row of $X$ is not zero,  there exist $a_{h}\in K$, $h=1, \ldots , u$, such that the coordinates of $Q$ verify the following conditions:
$$\label{rk1}\tilde{x}_{h;i,j}=a_{h}\tilde{x}_{1;i,j}\; \; \hbox{ and }\; \; x_{h,l}=a_{h}x_{1,l}$$ 
for $h=1, \ldots , {d-t+1\choose 2}$,  $i=1,2,3$, $j=1\ldots , t-k+1$ and $l=1, \ldots ,2k-t$.
\\
We are looking for a point $P:[z_{1}, \ldots , z_{u}]$ such that if the coordinates of $Q$ are as above, then $P$ and $Q$ verify the systems (\ref{Sh}). 
If $Q$ verifies (\ref{Sh}), then the coordinates of $P$ are such that:
$$\left(\begin{array}{cccc}
0&\cdots &\cdots&0\\
-a_{2}& a_{1}& \cdots &0\\
\vdots&    &\ddots &\\
-a_{u}& 0&\cdots& a_{1}
\end{array}\right)
\left( \begin{array}{c}
z_{1}\\
\vdots\\
z_{u}
\end{array}\right)=\underline{0},$$
that is to say $a_{h}z_{1}=z_{h}$ for $h=2 , \ldots , u$.
\\
The solution of such a system is the point $P$ we are looking for, i.e. $P:[a_{1}, \ldots , a_{u}]$.
\end{proof}

\subsection{The ideal of projections of Veronese surfaces from points}

\begin{theorem}\label{set} Let $X_{Z,d}$ be the projection of the Veronese $d$-uple embedding of $\PP 2$ from $Z=\{P_{1}, \ldots , P_{s}\}$ general points, $s\leq {d\choose 2}$. Then the equations (XM) and (Cat) together with either (\ref{equinx}) if $t/2\leq k \leq t$, or (\ref{E2}) if $0\leq k < t/2$, describe set theoretically $X_{Z,d}$.
\end{theorem}

\begin{proof}
Obviously $X_{Z,d}$ is contained in the support of the variety defined by the equations 
in statement of the theorem. 

In order to prove the other inclusion we need to prove that if a point $Q$ verifies all the equations required in the statement, then $Q\in X_{Z,d}$.

If $Q:[\tilde{x}_{h;i,j},x_{h,l}]$ annihilates the equations (XM), then, by Proposition \ref{PQ}, there exists a point $P:[z_{1}, \ldots , z_{u}]$ such that $P$ and $Q$ verify the systems (\ref{Sh}). Solving those systems in the variables $\tilde{x}_{h;i,j},x_{h,l}$ allows to write the point $Q$ depending on the $z_{1}, \ldots , z_{u}$. We do not write the computations for sake of simplicity, but  what it turns out is that there exist $\tilde{c}_{i,j}, c_{l}\in K$, with $i=1,2,3$, $j=1, \ldots ,t-k+1$ and $l=1, \ldots , 2k-t$ (only if $t/2 \leq k \leq t$) such that the coordinates $\tilde{x}_{h;i,j},x_{h,l}$ of $Q$ are  $\tilde{x}_{h;i,j}=\tilde c_{i,j}z_{h}$ and $x_{h,l}=c_{l}z_{h}$:
$$Q:[\tilde{x}_{h;i,j},x_{h,l}]=[\tilde c_{i,j}z_{h},c_{l}z_{h}].$$

Since such a $Q$, by hypothesis, verifies the equations (Cat), then there exists an unique point $R:[w_{1}, w_{2},w_{3}]\in \PP 2$ such that $z_{1}=w_{1}^{d-t-1}, z_{2}=w_{1}^{d-t-2}w_{2}, \ldots ,w_{3}^{d-t-1}$, therefore 
$$Q:[\tilde c_{i,j}\w^{\underline{\alpha}}, c_{l}\w^{\underline{\alpha}}]$$
with $|\underline{\alpha}|=d-t-1$.

Assume that  $R\notin Z$, that corresponds to assuming that $Q$ lies in the open set given by the image of $\varphi_{\tilde{J}_{d}}$ minus the exceptional divisors of $Bl_{Z}(\PP 2)$.

Now, if $t/2\leq k \leq t$, the point $Q$ verifies also the equations (\ref{equinx}), while if $0\leq k <t/2$ the point $Q$ verifies the equations (\ref{E2}). 
Therefore if $t/2\leq k \leq t$, then $\tilde c_{i,j}=b\tilde G_{i,j}$ and $c_{l}=bG_{l}$ for $i=1,2,3$, $j=1, \ldots ,t-k+1$ and $l=1, \ldots , 2k-t$;
if $0\leq k <t/2$, then $\tilde c_{i,j}=b\tilde G_{i,j}$ for $i=1,2,3$ and $j=1, \ldots ,t-k+1$, for some $b\in K$.
This proves that $Q\in X_{Z,d}$.
\end{proof}

Now we want to construct a weak generic hypermatrix of indeterminates $\mathcal{A}$ in the variables $\tilde {x}_{h;i,j},x_{h,l}$ in such a way that the vanishing of its $2$-minors coincide with the equations (XM) and (Cat). Then $I_{2}({\mathcal{A}})$ will be a prime ideal because of Proposition \ref{gen1.12}. so it will only remain to show that the generators of $I_{2}(\mathcal{A})$, together with the equations either (\ref{equinx}) or (\ref{E2}), are generators for the defining ideal of $X_{Z,d}$.
\\
\\
Let $C=(c_{i_{1},i_{2}})\in M_{3,d-t-3}(K)$ be the Catalecticant matrix defined in (\ref{C}). Let the $\tilde{x}_{h;i,j}$ and the $x_{h,l}$ be defined as in (\ref{x}). For all $i_{1}=1,2,3$, $i_{2}=1, \ldots , d-t-3$ and $i_{3}=1, \ldots , r$ where $r=2t-k+3$ if $t/2\leq k \leq t$ and $r=3(t-k+1)$ if $0\leq k <t$, construct the hypermatrix 
\begin{equation}\label{A}
\mathcal{A}=(a_{i_{1},i_{2},i_{3}})
\end{equation} 
in the following way: 
\begin{description}
\item[$a_{i_{1},i_{2},i_{3}}=\tilde{x}_{h,i,j}$] if $c_{i_{1},i_{2}}=z_{h}$ for $h=1, \ldots , {d-t+1\choose 2}$, and $i_{3}=1, \ldots , 3(t-k+1)$ is the position of the index $(i,j)$ after having ordered the $\tilde{G}_{i,j}$ with the lexicographic order,
\item[$a_{i_{1},i_{2},i_{3}}=x_{h,i_{3}-3(t-k+1)}$] if $c_{i_{1},i_{2}}=z_{h}$ for  $h=1, \ldots , {d-t+1\choose 2}$ and $i_{3}-3(t-k+1)=1, \ldots , 2k-t$ if $t/2 \leq k \leq t$. 
\end{description}

\begin{propos}\label{A weak} The hypermatrix $\mathcal{A}$ defined in (\ref{A}) is a weak generic hypermatrix of indeterminates.
\end{propos}

\begin{proof}  We need to verify that all the properties of weak generic hypermatrices hold for such an $\mathcal{A}$.
\begin{enumerate}
\item The fact that $\mathcal{A}=(\tilde{x}_{h;i,j}, x_{h,l})$ is a hypermatrix of indeterminates is obvious.
\item The variable $\tilde{x}_{1,1,1}$ appears only in position $a_{1,1,1}$.
\item The ideals of $2$-minors of the sections obtained fixing the third index of $\mathcal{A}$ are prime ideals because those sections are Catalecticant matrices and their $2$-minors are the equations of a Veronese embedding of $\PP 2$. The sections obtained fixing either the index $i_{1}$ or the index $i_{2}$ are generic matrices of indeterminates, hence their $2$-minors generate prime ideals.
\end{enumerate}
\end{proof}

\begin{corol}\label{primeness} Let $\mathcal{A}$ be defined as in (\ref{A}). The ideal $I_{2}({\mathcal{A}})$ is a prime ideal.
\end{corol}

\begin{proof}
This corollary is a consequence of Proposition \ref{A weak} and of Proposition \ref{gen1.12} .
\end{proof}

Now, we need to prove that the vanishing of the $2$-minors of the hypermatrix $\mathcal{A}$ defined in (\ref{A}) coincide with the equations (XM) and (Cat). 

\begin{theorem}\label{end surf} Let $X_{Z,d}$ be as in Theorem \ref{set}, then  the ideal $I(X_{Z,d})\subset K[\tilde{x}_{h;i,j}, x_{h,l}]$, with $h=1, \ldots , {d-t+1\choose 2}$,  $i=1,2,3$, $j=1\ldots , t-k+1$ and $l=1, \ldots ,2k-t$ is generated by all the $2$-minors of the hypermatrix $\mathcal{A}$ defined in (\ref{A}) and the linear formss appearing either in 
 (\ref{equinx}) if $t/2\leq k \leq t$ or in (\ref{E2}) if $0\leq k <t/2$.
\end{theorem}

\begin{proof}  In Corollary \ref{primeness} we have shown that $I_{2}({\CA})$ is a prime ideal; in Theorem \ref{set} we have proved that the equations (XM), (Cat) and either the equations (\ref{equinx}) if $t/2\leq k \leq t$ or the equations (\ref{E2}) if $0\leq k <t/2$ define $X_{Z,d}$ set-theoretically. Then we need to prove that the vanishing of the $2$-minors of $\mathcal{A}$ coincide with the equations (XM) and (Cat) and that either $(I_{2}({\CA}),(E_{1}))$ for $t/2 \leq k \leq t$, or $(I_{2}({\CA}),(E_{2}))$ is actually equal to $I(X_{Z,d})$  for $0\leq k \leq t/2$.

Denote with $I$ the ideal defined by $I_{2}({\mathcal{A}})$ and the polynomials appearing either in  (\ref{equinx}) in one case or in (\ref{E2}) in the other case. Denote also $\mathcal{V}$ the variety defined by $I$.

The inclusion ${\mathcal{V}}\subseteq X_{Z,d}$ is obvious because, by construction of $\mathcal{A}$, the ideal $I_{2}({\mathcal{A}})$ contains the equations (XM) and (Cat), therefore $I$ contains the ideal defined by $(XM)$, (Cat) and either (\ref{equinx})  or  (\ref{E2}).

For the other inclusion it is sufficient to verify that each $2$-minor of $\mathcal{A}$ appears either in (XM) or in (Cat). This is equivalent to prove that if $Q\in X_{Z,d}$ then $Q \in{\mathcal{V}}$, i.e. if $Q\in X_{Z,d}$ then $Q$ annihilates all the polynomials appearing in $I$.

An element of $I_{2}({\mathcal{A}})$ with ${\mathcal{A}}=(a_{i_{1},i_{2},i_{3}})$ is, by definition of a $2$-minor of a hypermatrix, one of the following:
\begin{enumerate}
\item $a_{i_{1},i_{2},i_{3}}a_{j_{1},j_{2},j_{3}}-a_{j_{1},i_{2},i_{3}}a_{i_{1},j_{2},j_{3}}$,
\item $a_{i_{1},i_{2},i_{3}}a_{j_{1},j_{2},j_{3}}-a_{i_{1},j_{2},i_{3}}a_{j_{1},i_{2},j_{3}}$,
\item $a_{i_{1},i_{2},i_{3}}a_{j_{1},j_{2},j_{3}}-a_{i_{1},i_{2},j_{3}}a_{j_{1},j_{2},i_{3}}$.
\end{enumerate}
We write for brevity $z_{i_{1},i_{2}}$ instead of $z_{h}$ if $(i_{1},i_{2})$ is the position occupied by $z_{h}$ in the catalecticant matrix $C$ defined in (\ref{C}). We also rename the $\tilde{G}_{i,j}$'s  and the $G_{l}$'s with $\overline{G}_{l}:=\tilde {G}_{i,j}$ if $l=1, \ldots , 3(t-k+1)$ is the position of $(i,j)$ ordered with the lexicographic order, and $\overline{G}_{l}:=G_{l -3(t-k+1)}$ if $l-3(t-k+1)=1, \ldots , 2k-t$.
\\
With this notation we evaluate those polynomials on $Q\in X_{Z,d}$.
\begin{enumerate}
\item $a_{i_{1},i_{2},i_{3}}a_{j_{1},j_{2},j_{3}}-a_{j_{1},i_{2},i_{3}}a_{i_{1},j_{2},j_{3}}=\overline{G}_{i_{3}}\overline{G}_{j_{3}}(z_{i_{1},i_{2}}z_{j_{1},j_{2}}-z_{j_{1},i_{2}}z_{i_{1},j_{2}})$ that vanishes on $X_{Z,d}$ because, by definition, $z_{1}=w_{1}^{d-t-1}$, $z_{2}=w_{1}^{d-t-2}w_{2}$, $\ldots$, $z_{u}=w_{3}^{d-t-1}$, hence the $z_{i,j}$'s vanish on the equations of the Veronese surface $Y_{2,d-t-1}$. The polynomial inside the parenthesis above is a minor of the catalecticant matrix defining such a surface, so the minor of $\mathcal{A}$ that we are studying vanishes on $X_{Z,d}$.
\item The above holds also for the case $a_{i_{1},i_{2},i_{3}}a_{j_{1},j_{2},j_{3}}-a_{i_{1},j_{2},i_{3}}a_{j_{1},i_{2},j_{3}}$.
\item $a_{i_{1},i_{2},i_{3}}a_{j_{1},j_{2},j_{3}}-a_{i_{1},i_{2},j_{3}}a_{j_{1},j_{2},i_{3}}=z_{i_{1},i_{1}}\overline{G}_{i_{3}}z_{j_{1},j_{2}}\overline{G}_{j_{3}}-z_{i_{1},i_{2}}\overline{G}_{j_{3}}z_{j_{1},j_{2}}\overline{G}_{i_{3}}=0$, evidently.
\end{enumerate}
This proves that the vanishing of the $2$-minors of $\CA$ coincides with the equations (XM) and (Cat).

For the remaining part of the proof, we work as in (\cite{Ha1}),  proof of Theorem 2.6.

Consider, with the previous notation, the sequence of surjective ring homomorphisms:
$$\begin{array}{ccccc}K[x_{i,j}]&\stackrel{\phi}{\rightarrow} &K[\w^{\underline{\alpha}}t_{j}]&\stackrel{\psi}{\rightarrow}&K[\w^{\underline{\alpha}}\overline{G}_{j}]\\
x_{i,j}&\mapsto & \w^{\underline{\alpha}}t_{j}&\mapsto&\w^{\underline{\alpha}}\overline{G}_{j}
\end{array}$$
where the exponent $\underline{\alpha}$ appearing in $\phi(x_{i,j})$ is the triple-index that is in position $i$ after having ordered the $\w$'s with the lexicographic order.

The ideal $I_{2}({\CA})$ is prime, so $I_{2}({\CA})\subseteq \ker(\phi )$.

Let $J\subset K[\w^{\underline{\alpha}}t_{j}]$ be the ideal generated by the images via $\phi$ of the equations appearing either in (\ref{equinx}) or in (\ref{E2}). The generators of $J$ are zero when $t_{j}=\overline{G}_{j}$, then $K[\w^{\underline{\alpha}}t_{j}]/J\simeq K[\w^{\underline{\alpha}}\overline{G}_{j}]$. Hence $J=\ker (\psi)$.
\\
Since it is almost obvious that a set of generators for $\ker(\psi \circ \phi)$ can be chosen as the generators of $\ker(\phi)$ together with the preimages via $\phi$ of the generators of $\ker(\psi)$, then $I= \ker (\psi \circ \phi)$. This is equivalent to the fact  that $I(X_{Z,d})=I$. 
\end{proof}

\section{Projection of Veronese varieties}\label{proiezioni veronese variet}

Here we want to generalize the results of the previous section to projections of Veronese varieties from a particular kind of irreducible and smooth varieties $V\subset \PP n$   of codimension $2$.

Since we want to generalize the case of $s$ general points in $\PP 2$, we choose $V$ of degree $s={t+1 \choose 2}+k \leq {d\choose 2}$ for some non negative integers $t$, $k$, $d$ such that $0<t<d-1$ and $0\leq k \leq t$.

Moreover we want to define the ideal $I(V)\subset K[x_{0}, \ldots , x_{n}]$ of $V$ as we defined  $J\subset K[x_{0}, x_{1},x_{2}]$ in Section \ref{4.1} (with the obvious difference that the elements of $I(V)$ belong to $K[x_{0}, \ldots , x_{n}]$ instead to $K[x_{0},x_{1},x_{2}]$). To be precise: let $L_{i,j}\in K[x_{0}, \ldots , x_{n}]_{1}$ be generic linear forms, and let $Q_{h,l}\in K[x_{0}, \ldots , x_{n}]_{2}$ be generic quadratic forms for $i,h=1, \ldots , k$, $j=1, \ldots , 2k-t$ and $l=1, \ldots , t-k+1$ if $t/2\leq k \leq t$; and for $i=1, \ldots , t-2k$, $j,l=1, \ldots , t-k+1$ and $h=1, \ldots , k$ if $0\leq k < t/2$. Define the matrix $\mathcal L$ either as in (\ref{Ltpiccolo}) or as in (\ref{Ltgrande}). The forms $F_{j}$ and $G_{l}$ are the maximal minors of $\mathcal L$ as previously. For each index $j$ there exist $n+1$ forms $\tilde G_{i,j}= w_{i}F_{j}$ with $i=0, \ldots ,n$, because now $\underline{w}=(w_{0}, \ldots , w_{n})$. Then the degree $d$ part of $I(V)$ is defined as $J_{d}$ in (\ref{Jd}) if $t/2 \leq k \leq t$ and as $J_{d}$ in (\ref{Jd2}) if $0\leq k < t/2$.
\\
This will be the scheme:
\begin{equation}\label{IV}
(V,I(V))\subset (\PP n, K[x_{0}, \ldots , x_{n}]).
\end{equation}
\\
\\
\Obs Let $W\subset \PP n$ be a variety of codimension $2$ in $\PP n$. Let $Y_{W}$ be the blow up of $\PP n$ along $W$. Let $E$ be the exceptional divisor of the blow up and $H$ the strict transform of a generic hyperplane. In \cite{Cop} (Theorem 1) it is proved that if $W$ is smooth, irreducible and scheme-theoretically generated in degree at most $\lambda\in \ZZ^{+}$, then $|dH-E|$ is very ample on the blow up $Y_{W}$ for all $d\geq \lambda +1$.
\\
\\
\Obs If $\deg(V)=s={t+1 \choose 2}+k \leq {d\choose 2}$, $0<t<d-1$ and $0\leq k \leq t$,
then $I(V)$  is generated in degrees $t$ and $t+1$.
\\
\\
A consequence of those remarks is the following:

\begin{propos}\label{bound b}
Let $V\subset \PP n$ be defined as in (\ref{IV}), and let $d > t+1$. If $E$ is the exceptional divisor of the blow up $Y_{V}$ of $\PP n$ along $V$  and $H$ is the strict transform of a generic hyperplane of $\PP n$, then $|dH-E|$ is very ample.
\end{propos}

Let $X_{V,d}\subset \PP {}(H^{0}({\mathcal O}_{Y_{V}}(dH-E)))$ be the image of the morphism associated to $|dH-E|$.

The arguments and the proofs used to study the ideal $I(X_{Z,d})$ in the previous section can all be generalized  to $I(X_{V,d})$ if $d>t+1$, $\deg (V) = {t+1 \choose 2}+k\leq {d\choose 2}$. \\
Now let $S'$ be the coordinate ring on $\PP {}(H^{0}({\mathcal O}_{Y_{V}}(dH-E)))$, constructed as $K[\tilde x_{i,j}, x_{h,l}]$ in the previous section: $S'=K[\tilde x_{i,j}, x_{h,l}]$ with $i=0, \ldots ,n$; $j=1, \ldots ,t-k+1$; $h=1, \ldots , {n+d-t-1 \choose 2}$ and $l=1, \ldots , 2k-t$ only if $t/2\leq k \leq t$ (in the other case the variables $x_{h,l}$ do not exist).

Let $(E')$ and $(E'')$ be the equations in $S'$ corresponding to  (\ref{equinx}) and (\ref{E2}), respectively.

Let $C'$ be the catalecticant matrix used to define the Veronese variety $Y_{n,d-t-1}$.

The hypermatrix $\CA '$ that we are going to use in this case is the obvious generalization of the hypermatrix $\CA $ defined in (\ref{A}); clearly one has to substitute $C$ with $C'$.

Now the proof of the fact that $I_{2}(\CA ')\subset S'$ is a prime ideal is analogous to that one of Corollary \ref{primeness}, and pass through the fact that $\CA '$ is a weak generic hypermatrix, hence we get the following:

\begin{theorem}
Let $(V,I(V))\subset (\PP n, K[x_{0}, \ldots ,x_{n}])$ be defined as in (\ref{IV}), let $Y_{V}$ be the blow up of $\PP n$ along $V$ and let $X_{V,d}$ be the image of $Y_{V}$ via $|dH-E|$, where $d> t+1$, $\deg (V)={t+1\choose 2}+k\leq {d\choose 2}$, $H$ is a generic hyperplane section of $\PP n$ and  $E$ is the exceptional divisor of the blow up. The ideal $I(X_{V,d})\subset S'$ is generated by all the $2$-minors of the hypermatrix $\CA '$ and the polynomials appearing either in ($E'$) if $t/2 \leq k \leq t$ or in ($E''$) if $0\leq k < t/2$, where $S'$, $\CA '$, ($E'$) and $(E'')$ are defined as above.
\end{theorem}


\begin{thebibliography}{Dillo 83}
\bibitem[\textbf{AR}]{AR} E.S. Allman, J. A. Rhodes, \emph{Phylogenetic invariants for the general Markov model of sequence mutation.} Math. Biosci. \textbf{186}  (2003), 113-144. MR2024609 (2004j:92048).
\bibitem[\textbf{BCS}]{BCS} P. B\"urgisser, M. Clausen, M. A. Shokrollahi, \emph{Algebraic complexity theory. With the collaboration of T. Lickteig.} Grundlehren der Mathematischen Wissenschaften [Fundamental Principles of Mathematical Sciences], 315. Springer-Verlag, Berlin, 1997. xxiv+618 pp. ISBN: 3-540-60582-7. MR1440179 (99c:68002).
\bibitem[\textbf{Bo}]{Bo} C. Bocci, \emph{Topics on Phylogenetic Algebraic Geometry}, Expositiones Mathematicae, \textbf{25}, no. 3 (2007), 235-259.
\bibitem[\textbf{Br}]{Br} R. Bro, \emph{PARAFAC, tutorial and applications}, Chemom. Intel. Lab. Syst., \textbf{38}, pp. 149Ð171, (1997). 
\bibitem[\textbf{BZ}]{BZ} I. Bengtsson, \.Zyczkowski, \emph{Geometry of quantum states. An introduction to quantum entanglement.} Cambridge University Press, Cambridge, (2006). xii+466 pp. ISBN: 978-0-521-81451-5; 0-521-81451-0. MR2230995 (2007k:81001).
\bibitem[\textbf{CHTV}]{CHTV} A. Conca, J. Herzog, N.V. Trung, G. Valla, \emph{Diagonal subalgebras of bigraded algebras and embeddings of blow-ups of projective spaces.} American Journal of mathematics, \textbf{119} (1997), 859-901. MR1465072 (99d:13001).
\bibitem[\textbf{CKP}]{CKP} J. D. Caroll, J. B. Kruskal, S. Pruzansky, \emph{Candelinc: A general approach to multidimensional analysis of many-way arrays with linear constraints on parameters}, Psychometrika, \textbf{45}, no. 1, pp. 3Ð24, Mar. (1980). 
\bibitem[\textbf{CGG1}]{CGG1} M. V. Catalisano, A. V. Geramita, A. Gimigliano \emph{Ranks of tensors, secant varieties of Segre varieties and fat points}. Linear Algebra Appl. \textbf{355} (2002), 263-285. MR1930149 (2003g:14070).
\bibitem[\textbf{CGG2}]{CGG2} M. V. Catalisano, A. V. Geramita, A. Gimigliano \emph{Higher Secant Varieties of Segre-Veronese varieties}. Atti del Convegno: Varieties with unexpected properties. Siena,  Giugno 2004. BERLIN: W. de Gruyter.  (2005 pp. 81 - 107)   MR2202248 (2007k:14109a).
\bibitem[\textbf{CGO}]{CGO} C. Ciliberto, A. V. Geramita, F. Orecchia \emph{Perfect varieties with defining equations of high degree.} Boll. U.M.I., Vol \textbf{1-B} no. 7,  (1987), 633-647. MR0916283 (88j:14061).
\bibitem[\textbf{Com}]{Com} P. Comon, \emph{Tensor decompositions: state of the art and applications.} Mathematics in signal processing, V (Coventry, 2000), 1--24, Inst. Math. Appl. Conf. Ser. New Ser., 71, Oxford Univ. Press, Oxford, 2002. MR1931400.
\bibitem[\textbf{Cop}]{Cop} M. Coppens \emph{Embedding of blowing-ups.} Sem. di Geometria 1991/93, Univ. di Bologna, Bologna (1994). MR1265754 (94m:14019).
\bibitem[\textbf{Ge}]{Ge} A. V. Geramita, \emph{Catalecticant varieties.} Commutative algebra and algebraic geometry (Ferrara), 143--156, Lecture Notes in Pure and Appl. Math., 206, Dekker, New York, 1999. 
\bibitem[\textbf{GG}]{GG} A. V. Geramita, A. Gimigliano \emph{Generators for the defining ideal of certain rational surfaces.} Duke Math. J. \textbf{62} (1991), 61-83. MR1104323 (92f:14031).
\bibitem[\textbf{Gi}]{Gi} A. Gimigliano \emph{On Veronesean surfaces.} Proc. Konin. Ned. Acad. van Wetenschappen, (A) \textbf{92} (1989), 71-85. MR0993680 (92a:14034).
\bibitem[\textbf{GL}]{GL} A. Gimigliano, A. Lorenzini \emph{On the Ideal of Veronese Surfaces.} Can. J. Math., Vol \textbf{45} no. 4,  (1993), 758-777. MR1227658 (94f:14031).
\bibitem[\textbf{Gr}]{Gr} R. Grone, \emph{Decomposable tensors as a quadratic variety.} Proc. of Amer. Math. \textbf{43} no. 2, (1977), 227-230. MR0472853 (57 $\sharp$12542).
\bibitem[\textbf{GSS}]{GSS} L. D. Garcia, M. Stillman, B. Strumfels, \emph{Algebraic Geometry of bayesian network.} J. Simbolic. Comp. \textbf{39}  (2005), 331-355. MR2168286 (2006g:68242).
\bibitem[\textbf{Ha1}]{Ha1} H. T. H\`a, \emph{Box-shaped matrices and the defining ideal of certain blowup surface.} Journal of Pure and Applied Algebra. \textbf{167}  no. 2-3, (2002), 203-224. MR1874542 (2002h:13020).
\bibitem[\textbf{Ha2}]{Ha2} H. T. H\`a, \emph{On the Rees algebra of certain codimension two perfect ideals.} Manuscripta Matematica. \textbf{107}, (2002), 479-501. MR1906772 (2003d:13002).
\bibitem[\textbf{HR}]{HR} S. Ho\c sten, S. Ruffa, \emph{Introductory notes to algebraic statistics.} Rend. Istit. Mat. Univ. Trieste. \textbf{37}  no. 1-2, (2005), 39-70. MR2227048.
\bibitem[\textbf{Lak}]{Lak} J. A. Lake, \emph{A rate-independent technique for analysis of nucleic acid sequences: evolutionary parsimony.} Mol. Biol. Evol. \textbf{4}  no. 2, (1987), 167-191. 
\bibitem[\textbf{Lan}]{Lan} J. M. Landsberg, \emph{Geometry and the complexity of Matrix Multiplication}, preprint: arXiv:cs/0703059v1 [cs.CC].
\bibitem[\textbf{Li}]{Li} T. Lickteig, \emph{Typical tensor rank.} Linear Algebra Appl. \textbf{69}, (1985), 95-120. MR0798367 (87f:15017).
\bibitem[\textbf{LM}]{LM} J. M. Landsberg, L. Manivel \emph{On the ideals of secant varieties of Segre varieties}. Found Comput. Math. \textbf{4} (2004), no. 4, 397-422. MR2097214 (2005m:14101).
\bibitem[\textbf{Lo}]{Lo} A. Lorenzini, \emph{Betti numbers of perfect homogeneous ideals.} J. of Pure and Applied Alg. \textbf{60}, (1989), 273-288. MR1021852 (90i:13022).
\bibitem[\textbf{LW}]{LW} J. M. Landsberg, J. Weyman \emph{On the ideals and singularities of secant varieties of Segre varieties}. preprint math.AG/0601452.
\bibitem[\textbf{MU}]{MU} S. Morey, B. Ulrich, \emph{Rees algebra of ideals with low codimension.} Proceedings of AMS. \textbf{124}, (1996), 3653-3661. MR1343713 (97e:13006).
\bibitem[\textbf{AOP}]{AOP} H. Abo, G. Ottaviani, P. Peterson, \emph{Induction for secant varieties of Segre varieties.} preprint math.AG/0607191. 
\bibitem[\textbf{Pa}]{Pa} A. Parolin \emph{Variet\`a Secanti alle Variet\`a di Segre e di Veronese e Loro Applicazioni}, Tesi di dottorato, Universit\`a di Bologna, A.A. 2003/2004.
\bibitem[\textbf{PS}]{PS} L. Pachter, B. Strumfels, \emph{Algebraic statistics for computational biology.} Cambridge University Press, New York (2005), MR MR2205865 (2006i:92002). 
\bibitem[\textbf{Pu}]{Pu} M. Pucci, \emph{The Veronese variety and Catalecticant matrices.} J. Algebra. \textbf{202} no. 1, (1998), 72-95. MR1614174 (2000c:14071).
\bibitem[\textbf{Sh}]{Sh} D.W. Sharpe, \emph{On certain polynomial ideals defined by matrices.} Quart. J. Math.Oxford \textbf{15} no. 2, (1964), 155-175. MR0163927 (29 $\sharp$ 1226).
\bibitem[\textbf{SS}]{SS} J.G. Semple; L. Roth, \emph{Introduction to Algebraic Geometry.} Oxford, at the Clarendon Press, 1949. xvi+446 pp. MR0034048 (11,535d).
\bibitem[\textbf{St}]{St} V. Strassen, \emph{Relative bilinear complexity and matrix multiplication.} J.  Reine Angew. Math. \textbf{375/376}, (1987), 406-443. MR0882307 (88h:11026).
\bibitem[\textbf{Wa}]{Wa} K. Wakeford, \emph{On Canonical Forms.} Proc. London Math. Soc. \textbf{18} (1918-19), 403-410.
\bibitem[\textbf{Za}]{Za} F. L. Zak, \emph{Tangents and secant of algebraic varieties}, Translation of Mathematical Monographs, \textbf{127}, American Mathematical Society, Providence, RI, 1993, Translated from the Russian manuscript by the author. MR1234494 (94i:14053).
\end{thebibliography}
\end{document}